\documentclass[final]{siamltex}
\usepackage{mathrsfs}
\usepackage{amsmath}
\usepackage{amsfonts}
\usepackage{amssymb}
\usepackage{graphicx}
\usepackage{geometry}
\usepackage{placeins}
\usepackage{mathrsfs}
\usepackage{color}
\usepackage{float}
\usepackage{subfigure}
\usepackage{threeparttable}

\newtheorem{rem}{Remark}[section]

\begin{document}

\title{New Ensemble Domain Decomposition Method for the Steady-state Random Stokes-Darcy Coupled Problems with Uncertain Parameters}
\author{
Chunchi Liu \thanks{%
School of Mathematical Sciences, East China Normal University,
Shanghai. \texttt{liuchunchimath@163.com}.}
\and Yao Rong \thanks{%
	College of Engineering Physics, Shenzhen Technology University,
	Shenzhen. \texttt{yr926@outlook.com}.}
\and Yizhong Sun \thanks{%
School of Mathematical Sciences, East China Normal University,
Shanghai. \texttt{bill950204@126.com}.}
\and Jiaping Yu \thanks{%
College of Science, Donghua University, Shanghai. \texttt{jpyu@dhu.edu.cn}.}
\and Haibiao Zheng \thanks{%
School of Mathematical Sciences, East China Normal University,
Shanghai Key Laboratory of Pure Mathematics and Mathematical
Practice, Shanghai. \texttt{hbzheng@math.ecnu.edu.cn}.
Partially
supported by Science and Technology
Commission of Shanghai Municipality (Grant Nos. 22JC1400900, 21JC1402500, 22DZ2229014). }}

 \maketitle

\begin{abstract}
 This paper presents two novel ensemble domain decomposition methods for fast-solving the Stokes-Darcy coupled models with random hydraulic conductivity and body force. To address such random systems, we employ the Monte Carlo (MC) method to generate a set of independent and identically distributed deterministic model samples. To facilitate the fast calculation of these samples, we adroitly integrate the ensemble idea with the domain decomposition method (DDM). This approach not only allows multiple linear problems to share a standard coefficient matrix but also enables easy-to-use and convenient parallel computing. By selecting appropriate Robin parameters, we rigorously prove that the proposed algorithm has mesh-dependent and mesh-independent convergence rates. For cases that require mesh-independent convergence, we additionally provide optimized Robin parameters to achieve optimal convergence rates. We further adopt the multi-level Monte Carlo (MLMC) method to significantly lower the computational cost in the probability space, as the number of samples drops quickly when the mesh becomes finer. Building on our findings, we propose two novel algorithms: MC ensemble DDM and MLMC ensemble DDM, specifically for random models. Furthermore, we strictly give the optimal convergence order for both algorithms. Finally, we present several sets of numerical experiments to showcase the efficiency of our algorithm.
\end{abstract}

\section{Introduction}

Simulating the coupling problem of free fluid flow and porous medium flow is crucial in engineering and geological applications. To model this multi-domain, multi-physics coupled systems, researchers have constructed various mathematical models, such as the Stokes-Darcy models \cite{DM00, DM000, LMLayton12, CGHW02, Cao03}, the Stokes-Brinkman model \cite{Xu09, Williamson10, Williamson11}, the Stokes-Darcy-transport/heat model \cite{Zhang01}, dual-porosity-Stokes model \cite{Hou06}, among others. Additionally, many effective methods have emerged to address these models, including Lagrange multiplier methods \cite{LMLayton12}, domain decomposition methods \cite{Sun04,Sun05,Discacciati18,Chen19,Cao20,Cao21,He22,Liu34}, and multi-grid methods \cite{Mu23,Cai24}, to name a few.


The domain decomposition method (DDM) can naturally decouple multi-domain multi-physics problems into smaller subproblems under appropriate interface conditions. It is easy to operate and facilitates parallel computing, so it has received extensive attention. Chen et al. \cite{Chen19} proposed a new parallel Robin domain decomposition method for the Stokes-Darcy model with Beavers-Joseph-Saffman (BJS) interface conditions. They demonstrated that the proposed algorithm has a mesh-independent geometric convergence rate with suitable choices of Robin parameters. However, they only provided a range to choose Robin parameters. Osmolovskii and Samarskii \cite{osm01, osm02} utilized the optimized Schwarz method to derive the optimized Robin parameters and achieve the optimal convergence rate of DDM for the Stokes-Darcy model with BJS interface conditions.


For effective simulations, parallel computing with DDM for accurate simulations is usually not feasible due to the fact that it is physically impossible to know the exact parameter values \cite{JiangN29}. In particular, the hydraulic conductivity tensor $\mathbb{K}$ is an uncertain parameter due to the complexity of geological structures in the porous media region, the large-scale real domain, and the small-scale natural randomness. One technique to deal with uncertainty is to model the uncertain parameters as a random function determined by a basic random field with a specified stipulated covariance structure (usually determined experimentally). This means we need to solve random partial differential equations (PDEs) \cite{FengXB}, and then such random PDEs can be solved using the Monte Carlo (MC) method. However, due to the slow convergence speed of the MC method, this method needs to consume a huge amount of computing resources to obtain accurate results. Considering such computational challenges, a class of ensemble methods has been developed \cite{ JiangN29, JiangNH27}. For deterministic models with different physical parameters or external forces, the above ensemble idea can enable multiple linear systems to share a common coefficient matrix, further improving computational efficiency. Ensemble algorithms have also been recently extended in \cite{JiangN30, JiangN31, SunF}. To overcome the high computational cost of MC in the probability space, Luo et al. \cite{Mul01} developed a MLMC ensemble scheme for solving parabolic equations with random coefficients. Recently, Yang et al. \cite{Mul02} proposed the multi-grid MLMC method for the Stokes-Darcy interface model with random hydraulic conductivity.

In this work, we mainly discuss an efficient ensemble domain decomposition method for the steady-state random Stokes-Darcy coupled model with uncertain hydraulic conductivity tensors. Different from \cite{Chen19}, the model we discuss is the random PDEs, which increase computational complexity and pose some analytical challenges. Inspired by \cite{JiangN29}, we follow the idea of the ensemble algorithm to overcome the difficulties caused by the randomness of parameters. It is worth noting that the authors of \cite{JiangN29} utilize the Gronwall lemma in their theoretical analysis of the unsteady coupled model, which does not apply to the steady-state model. Therefore, we choose the iterative technique to decouple the steady-state coupled system and present a necessary and important lemma (Lemma 4.1) for theoretical analysis. More specifically, we rigorously prove that our new algorithm converges in the continuous system while the Robin parameters satisfy $\gamma_{f} \leq \gamma_{p}$. An interesting result is that for $\gamma_{f} < \gamma_{p}$, the ensemble DDM has a mesh-independent geometric convergence rate. Furthermore, we note that the convergence of our algorithm is closely related to the choice of Robin parameters. To obtain the optimal convergence rate, we briefly deduce and analyze the optimized Robin parameters, which can also be proved to accelerate the convergence in the numerical experiments. Moreover, in the finite element approximation, when $\gamma_{f} = \gamma_{p}$, the convergence of the corresponding mixed finite element form for the ensemble DDM is demonstrated to be mesh-dependent. The classical MC method needs a larger number of samples in probability space and a smaller mesh size in physical space for a more accurate numerical approximation. This leads to a significant increase in the computational cost, which is the product of the number of samples and the cost of each sample. Therefore, we further adopt the multi-level Monte Carlo (MLMC) method, which cuts down the computational cost in the probability space dramatically, as the number of samples reduces fast while the mesh size decreases. To the best of our knowledge, this is the first work to concurrently consider reducing the number of samples and the cost of each sample, resulting in MLMC ensemble DDM. In the probability space, we utilize the MLMC method to reduce the number of samples. In the physical space, the algorithm integrates the ensemble idea and domain decomposition method to save the cost of each sample.

The structure of this paper is summarized as follows. Section 2 introduces a steady-state random Stokes-Darcy model with BJS interface conditions. In Section 3, we design a new parallel ensemble domain decomposition method with two Robin-type condition sets. In Section 4, we analyze the convergence of the proposed algorithm in the continuous system when the Robin parameters satisfy $\gamma_{f} \leq \gamma_{p}$. In this Section, we also obtain the mesh-independent convergence rate and optimized Robin parameters while $\gamma_{f} < \gamma_{p}$. The corresponding mixed finite element form for the ensemble DDM is proved to be mesh-dependent convergent for $\gamma_{f} = \gamma_{p}$ in Section 5. A brief introduction to the multil-evel Monte Carlo method is shown in Section 6. Section 7 presents some numerical tests to support our theoretical results.

\section{The Stokes-Darcy Model with BJS Interface Conditions}
\begin{figure}[htbp]
\centering
\includegraphics[width=90mm,height=45mm]{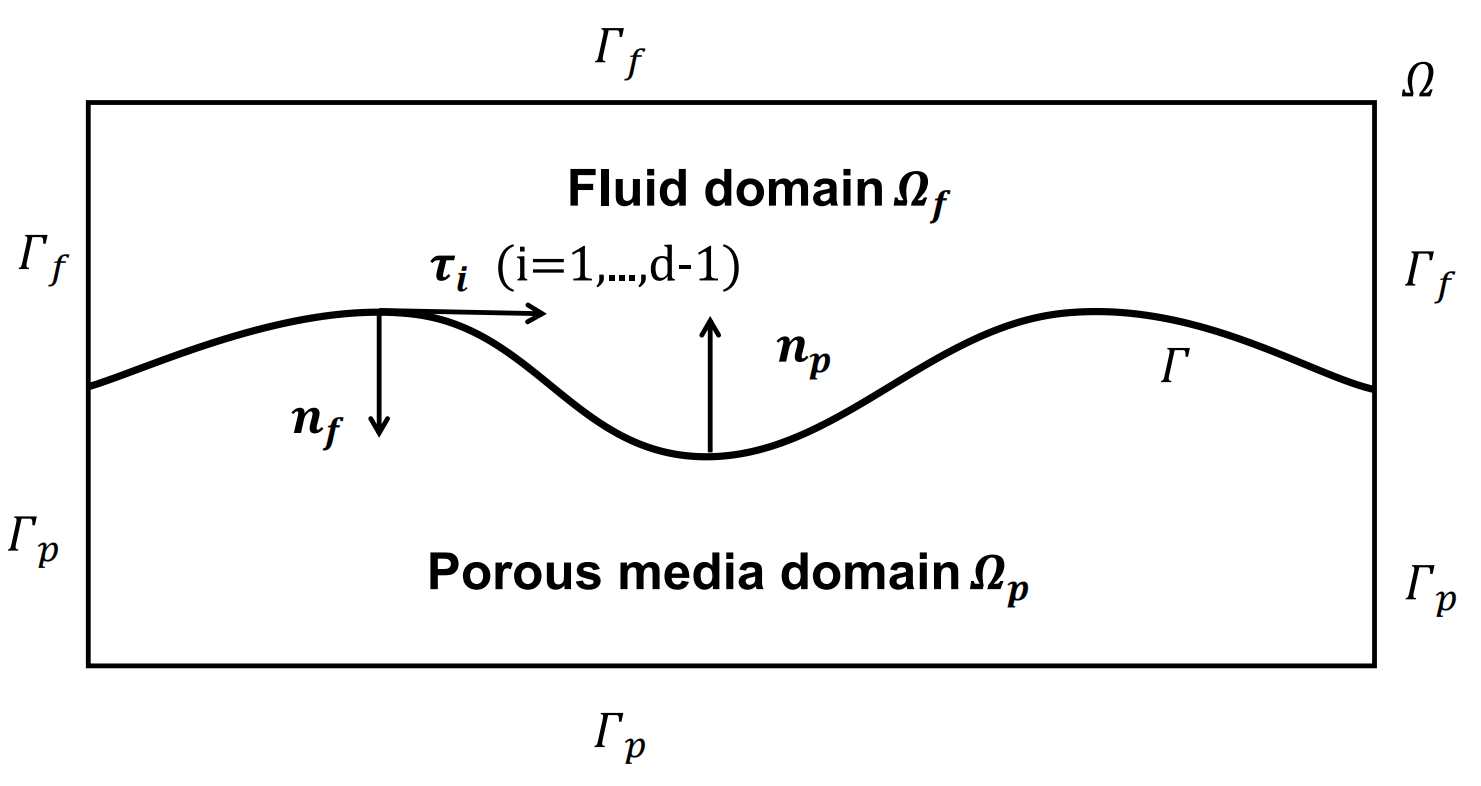}
\caption{The global domain $\overline{\Omega}$ consisting of the fluid domain
$\Omega_f$ and the porous media region $\Omega_p$ separated by the
interface $\Gamma$.}\label{dd}
\end{figure}

Consider two bounded non-overlapping domains $\Omega _{f},~\Omega_{p}\subset R^{d}$~$(d=2~\mathrm{or}~3)$ with an interface $\Gamma$, which denoted by $\Omega _{f}$ for fluid flow and $\Omega _{p}$ for porous media flow, namely $\Omega _{f}\cap \Omega _{p}={\emptyset }$, $\overline{\Omega }_{f}\cap\overline{\Omega }_{p}=\Gamma $. We define $\mathbf{n}_{f}$ and $\mathbf{n}_{p}$ are the unit outward normal vectors on $\partial \Omega _{f}$ and $\partial \Omega _{p}$, respectively. Those two symbols satisfy the equality relationship: $\mathbf{n}_{f}=-\mathbf{n}_{p}~\mathrm{on}~\Gamma $. Besides, the unit tangential vectors on the interface $\Gamma $ are represented by $\mathbf{\tau }_{i},~i=1, \cdots, d-1$, and $\Gamma_{f}=\partial \Omega _{f}\setminus \Gamma ,~\Gamma _{p}=\partial \Omega_{p}\setminus \Gamma $, see Fig. \ref{dd} for a sketch.

The fluid flow in domain $\Omega_{f} $ can be governed by the Stokes equations, finding the fluid velocity $\mathbf{u}_{f}$ and kinematic pressure $p_{f}$, such that
\begin{eqnarray}
 -\nabla \cdot \mathbb{T}(\mathbf{u}_{f},p_f)&=&\mathbf{f}\ \ \ \ \mathrm{in}~\Omega _{f},
\label{Stokes01} \\
\nabla \cdot \mathbf{u}_{f} &=&0\ \ \ \ \mathrm{in}~\Omega _{f},\label{Stokes02}
\end{eqnarray}
where $\mathbb{T}(\mathbf{u}_{f},p_f)=2\nu\mathbb{D}(\mathbf{u}_{f})-p_f\mathbb{I}$ is the stress tensor, $\mathbb{I}$ is identity matrix, $\mathbb{D}( \mathbf{u}_{f})=\frac{1}{2}(\nabla \mathbf{u}_{f}+(\nabla \mathbf{u}_{f})^{T})$ is the deformation tensor, $\nu $ means the kinematic viscosity of the fluid flow and $\mathbf{f}$ is the external body force.

In the porous media domain $\Omega_{p} $, the porous media flow
is governed by the following Darcy equation with the piezometric head $\phi_{p}$
\begin{eqnarray}
- \nabla \cdot (\mathbb{K}(\mathbf{x}) \nabla \phi_{p}) &=& 0\ \ \hspace{3mm} %
\mathrm{in}~\Omega _{p}, \label{Darcy01}
\end{eqnarray}
where $\mathbb{K}(\mathbf{x}) $ is the hydraulic conductivity tensor, which is physically impossible to determine its parameter values.
The piezometric head is $\phi_{p}=z+\frac{p_p}{\rho g}$, where $z$ is the height, $p_p$ indicates the dynamic pressure, $g$ denotes the gravitational acceleration, and $\rho$ respects density. 

We assume the fluid velocity $\mathbf{u}_{f} $ and the piezometric head $\phi_{p}$ satisfy the homogeneous Dirichlet boundary condition except on $\Gamma$, i.e., $\mathbf{u}_{f}=0~~\mathrm{on}~\Gamma _{f} $ and $\phi_{p}=0~~\mathrm{on}~\Gamma _{p} $. As for the interface $\Gamma $, the coupling conditions require to satisfy the conservation of mass, the balance of forces, and the tangential conditions on the fluid domain's velocity. Herein, we impose the Beavers-Joseph-Saffman (BJS) as the tangential condition on the interface $\Gamma $. The interface coupling conditions in this paper are shown as follows
\begin{eqnarray}
\mathbf{u}_{f}\cdot \mathbf{n}_{f}-\mathbb{K} \nabla \phi_{p}\cdot
\mathbf{n}_{p}
&=&0\ \ \ \ \ \ \ \ \ \ \ \ \ \ \hspace{7.5mm} \mathrm{on}~\Gamma ,  \label{interface1} \\
-\mathbf{n}_{f}\cdot (\mathbb{T}(\mathbf{u}_{f},p_f)\cdot \mathbf{n}_{f})
&=&g\phi_{p} \ \ \ \ \ \ \ \ \ \ \ \ \ \ \ \ \  \mathrm{on}~\Gamma , \label{interface2}\\
-\tau_{i}\cdot(\mathbb{T}(\mathbf{u}_{f},p_f)\cdot \mathbf{n}_{f})
&=&\frac{\alpha}{\sqrt{\tau_{i}\cdot \mathbb{K}\tau_{i}}} \tau _{i} \cdot \mathbf{u}_{f}
\hspace{3mm}  1\leq i \leq d-1\hspace{3mm} \mathrm{on}~\Gamma,\label{interface3}
\end{eqnarray}%
where $\alpha$ represents an experimentally determined positive parameter depending on the Darcy properties. 

Due to the hydraulic conductivity tensor $\mathbb{K}$ having small perturbations, it's not feasible to obtain the exact hydraulic conductivity values via limited equipment. Then we should discuss the Stokes-Darcy coupled model with the random hydraulic conductivity tensor $\mathbb{K}(\mathbf{x},\omega)$. Let $(X,\mathcal{F}, \mathcal{P})$ be a complete probability space, where $X$ is the set of outcomes, $\mathcal{F} \in 2^{X}$ is the $\sigma$-algebra of events and $\mathcal{P} : \mathcal{F} \rightarrow [0,1]$ is a probability measure.
The random Stokes-Darcy model could be written as: find the functions $\mathbf{u}_{f} \ : \ \Omega_f \times X \rightarrow \mathbb{R}^d \ (d=2,3)$,  ${p}_{f} \ : \ \Omega_f \times X \rightarrow \mathbb{R} $, and  ${\phi}_{p} \ : \ \Omega_p \times X \rightarrow \mathbb{R} $, such that it holds $\mathcal{P}-a.e.$ in $X$
\begin{eqnarray}
	-\nabla \cdot \mathbb{T}(\mathbf{u}_{f}(\mathbf{x},\omega),p_f(\mathbf{x},\omega))&=&\mathbf{f}(\mathbf{x},\omega)\ \ \ \ \ \ \ \ \mathrm{in}~\Omega_f \times X,
	\label{Stokes001} \\
	\nabla \cdot \mathbf{u}_{f}(\mathbf{x},\omega) &=&0\ \ \ \  \ \ \ \ \ \ \ \  \ \ \ \ \mathrm{in}~\Omega_f \times X,\label{Stokes002}\\
	\nabla \cdot (-\mathbb{K}(\mathbf{x},\omega) \nabla \phi_{p}(\mathbf{x},\omega))&=&0\ \ \ \ \ \ \  \ \ \ \ \ \ \ \ \ 
	\mathrm{in}~\Omega_p \times X ,  \label{Darcy001} 
\end{eqnarray}
where $\mathbf{f}(\mathbf{x},\omega)\ : \ \Omega_f \times X \rightarrow \mathbb{R}^d$, and $\mathbb{K}(\mathbf{x},\omega) $ is a continuous and bounded random function. For the random PDEs, the classical approach to solve this random system is the Monte Carlo method, and some details can be seen in \cite{JiangN29}. Hence, we assume that there are $J$ samples for the random $\mathbb{K}(\mathbf{x},\omega_j),~\mathbf{f}(\mathbf{x},\omega_j),~j=1,\cdots,J$ and the samples are independently and identically distributed. For simplicity, regarding $\mathbb{K}(\mathbf{x},\omega_j)$ as $\mathbb{K}_j(\mathbf{x})$, similarly denoting force $\mathbf{f}(\mathbf{x},\omega_j)$ as $\mathbf{f}_j(\mathbf{x})$, we have an ensemble of $J$ Stokes-Darcy system corresponding to $J$ different determined parameter sets $(\mathbf{f}_j,\mathbb{K}_j), j=1,\dots,J$ to be computed as 
\begin{eqnarray}
	-\nabla \cdot \mathbb{T}(\mathbf{u}_{j},p_j)&=&\mathbf{f}_{j}\ \ \ \ \ \ \ \ \ \ \ \ \ \ \ \mathrm{in}~\Omega _{f},
	\label{Stokes1} \\
	\nabla \cdot \mathbf{u}_{j} &=&0\ \ \ \  \ \ \ \ \ \ \ \  \ \ \ \ \mathrm{in}~\Omega _{f},\label{Stokes2}\\
	\nabla \cdot (-\mathbb{K}_j \nabla \phi_{j})&=&0\ \ \ \ \ \ \  \ \ \ \ \ \ \ \ \ 
	\mathrm{in}~\Omega_{p} .  \label{Darcy1} 
\end{eqnarray}

To solve the above problems, some Sobolev spaces are utilized as follows
\begin{eqnarray*}
\mathbf{X}_{f}&:=&\Bigl\{\mathbf{v}\in H^{1}(\Omega_{f})^{d}:\mathbf{v}=0~\mathrm{on}~\Gamma _{f}\Bigr\},\\
Q_{f}&:=&L^{2}(\Omega_{f}), \\
X_{p}&:=&\Bigl\{\psi\in H^1(\Omega_{p}): \psi =0~%
\mathrm{on}~\Gamma_{p}\Bigr\}.
\end{eqnarray*}
We denote the inner products for each domain by $(\cdot, \cdot )_{\Omega_{f}}$ and $(\cdot, \cdot )_{\Omega_{p}}$ respectively, and the corresponding $L^2$-norms by $||\cdot ||_{\Omega_{f}}$ and $||\cdot ||_{\Omega_{p}}$. Meanwhile, $\langle \cdot, \cdot \rangle$ is defined as the $L^2$ inner product on the interface $\Gamma$, with the $L^{2}(\Gamma )$ norm $||\cdot ||_{\Gamma }$. Then, the weak formulation of the Stokes-Darcy system is derived similarly to the equation (1.4)-(1.5) of the paper \cite{Chen19}.

\section{ Parallel Ensemble Domain Decomposition Method}
In this section, we mainly propose the ensemble domain decomposition method to decouple the Stokes-Darcy model. First, We discuss the Robin-type conditions for the free fluid system of the random Stokes equations. For given Robin parameter of Stokes $\gamma_f>0$, the corresponding function $\delta_{fj}$ is defined on $\Gamma$ 
\begin{eqnarray}
\delta_{fj}=\mathbf{n}_{f}\cdot(\mathbb{T}(\mathbf{u}_{j},p_{j})\cdot\mathbf{n}_f)+\gamma_{f}\mathbf{u}_{j}\cdot\mathbf{n}_{f}.\label{Robinstok1}
\end{eqnarray}
Then, the corresponding weak formulations for the Stokes model is: for $\delta_{fj}\in L^{2}(\Gamma )$, find $\mathbf{u}_{j}\in \mathbf{X}_{f}$ and $p_{j}\in Q_{f}$ satisfying the compatibility conditions (\ref{Robinstok1}) on the interface $\Gamma$, such that for all $(\mathbf{v},q)\in (\mathbf{X}_{f}, Q_{f})$
\begin{eqnarray}\label{decoupled-111}
	&&a_{f}(\mathbf{u}_{j},\mathbf{v})+b_f(p_{j},\mathbf{v})+\gamma_{f}\langle\mathbf{u}_{j}\cdot\mathbf{n}_{f},\mathbf{v}\cdot\mathbf{n}_f\rangle
	+\sum_{i=1}^{d-1}\langle\overline{\eta}_{i}\mathbf{u}_{j}\cdot\tau_{i},\mathbf{v}\cdot\tau_{i}\rangle\nonumber\\
	&&\hspace{3mm}=(\mathbf{f}_j,\mathbf{v})_{\Omega_f}-\sum_{i=1}^{d-1}\langle(\eta_{i j}-\overline{\eta}_{i})\mathbf{u}_{j}\cdot\tau_{i},\mathbf{v}\cdot\tau_{i}\rangle+\langle\delta_{fj},\mathbf{v}\cdot\mathbf{n}_{f}\rangle, \\
	&&b_f(\mathbf{u}_j,q)=0,\label{decoupled-222}
\end{eqnarray}
where the bilinear forms
\begin{eqnarray*}
	a_f(\mathbf{u}_{j},\mathbf{v})=2\nu(\mathbb{D}(\mathbf{u}_j),\mathbb{D}(\mathbf{v}))_{\Omega_f}, \hspace*{15mm}	b_f(\mathbf{v},q)=-(\triangledown\cdot \mathbf{v},q)_{\Omega_f}.
\end{eqnarray*}
Besides, some routine notations in the formulations (\ref{decoupled-111}) are denoted as
\begin{eqnarray*}
\eta_{ij}=\frac{\alpha}{\sqrt{\tau_{i}\cdot\mathbb{K}_j\tau_{i}}},\qquad\overline{\eta}_{i}=\frac{1}{J}\sum_{j=1}^{J}\eta_{i j}.
\end{eqnarray*}
For given Robin parameter of Darcy $\gamma_p>0$, the function $\delta_{pj}$ for the Darcy equation satisfies
\begin{eqnarray}
	\delta_{pj}=\gamma_{p}\mathbb{K}_{j}\triangledown\phi_{j}\cdot\mathbf{n}_{p}+g\phi_{j}.\label{Robindar1} 
\end{eqnarray}
Then, the weak formulation for the Darcy system is given: for $\delta_{pj}\in L^{2}(\Gamma )$, find $\phi_{j}\in X_{p}$ satisfying the compatibility conditions (\ref{Robindar1}) on the interface $\Gamma$, such that for all $\psi\in X_{p}$
\begin{eqnarray}
	\gamma_{p}(\overline{\mathbb{K}}\triangledown\phi_{j},\triangledown\psi)_{\Omega_{p}}+g\langle \phi_{j},\psi\rangle
	&=&\langle\delta_{pj},\psi\rangle-\gamma_{p}((\mathbb{K}_{j}-\overline{\mathbb{K}})\triangledown\phi_{j},\triangledown\psi)_{\Omega_{p}},\label{decoupled-333}
\end{eqnarray}
where
\begin{eqnarray*}
	\overline{\mathbb{K}}=\frac{1}{J}\sum_{j=1}^{J}\mathbb{K}_{j}.
\end{eqnarray*}

Let $||\cdot ||_2$ denote the $2$-norm of either vectors or matrices. Let $k_{j \min}(\mathbf{x}),\overline{k}_{\min}(\mathbf{x})$ be the minimum eigenvalue of the hydraulic conductivity tensor $\mathbb{K}_{j},\overline{\mathbb{K}}$ respectively, and $k_{j \max}(\mathbf{x})$ is the maximum eigenvalue of $\mathbb{K}_{j}$. In particular, we assume $\mathbb{K}_{j}$ and $\overline{\mathbb{K}}$ are symmetric. For simplicity, the spectral radius of the fluctuation of hydraulic conductivity tensor can be expressed as $||\mathbb{K}_{j}-\overline{\mathbb{K}}||_2=\rho_{j}^{\prime}(\mathbf{x})$.
We also assume $0\leq k_{j\min }\leq \lambda (\mathbb{K}_j)\leq k_{j\max}<\infty$.

\begin{lemma}
The interface conditions (\ref{interface1})-(\ref{interface3}) are equivalent to the Robin-type conditions (\ref{Robinstok1}),(\ref{Robindar1}) if and
only if $\delta_{fj}$ and $\delta_{pj}$ satisfy the following two compatibility conditions
on the interface $\Gamma$
\begin{eqnarray}
&&\delta_{pj}=\gamma_{p}\mathbf{u}_{j}\cdot\mathbf{n}_{f}+g\phi_{j}, \label{comp1}\\
&&\delta_{fj}=\gamma_{f}\mathbf{u}_{j}\cdot\mathbf{n}_{f}-g\phi_{j}.\label{comp2}
\end{eqnarray}
\end{lemma}

Due to the influence of random parameters, it is necessary to solve a great deal of algebraic equations with varying stiffness matrices. Moreover, to better estimate the uncertainty and sensitivity of the solution, we shall select more samples, which will essentially increase the computational cost. To improve the computational efficiency, we hope to get an algebraic structure similar to that shown below
\begin{eqnarray*}
	\begin{aligned}
		&A\left[\begin{array}{l l l}
			\mathbf{u}_{1} \\
			p_{1}
		\end{array}
		|\cdots| 
		\begin{array}{l}
			\mathbf{u}_{J} \\
			p_{J}
		\end{array}\right]=\left[RHS_{1}|\cdots| RHS_{J}\right], \hspace{4mm}
		&B\left[
			\phi_{1} | \cdots | \phi_{J} \right]=\left[RHS^*_{1}|\cdots| RHS^*_{J}\right],
	\end{aligned}
\end{eqnarray*}
which shares the same coefficient matrix. Such ensemble idea can make the coefficient matrices $A$ and $B$ only to use once efficient iterative solvers or direct solvers for fast computation. Hence, we propose the following parallel ensemble domain decomposition method.
\newline

\textbf{\emph{{Ensemble DDM Algorithm}}}

1. Initial values of $\delta_{pj}^0$, $\delta_{fj}^0$, $\mathbf{u}_j^0$ and $\phi_j^0$ are guessed, which can be zero.

2. For $n=0,1,2,\cdots,$ independently solve the ensemble Stokes and Darcy equations, such that for all $  (\mathbf{v},q;\psi)\in (\mathbf{X}_{f},Q_f;X_p)$, we obtain  $(\mathbf{u}_{j}^{n+1},p_{j}^{n+1};\phi_j^{n+1})\in  (\mathbf{X}_{f}, Q_{f};X_{p})$ by solving
\begin{eqnarray}\label{decoupled-1}
		&&a_{f}(\mathbf{u}_{j}^{n+1},\mathbf{v})+b_f(p_{j}^{n+1},\mathbf{v})
	+\gamma_{f}\langle\mathbf{u}_{j}^{n+1}\cdot\mathbf{n}_{f},\mathbf{v}\cdot\mathbf{n}_f\rangle
	+\sum_{i=1}^{d-1}\langle\overline{\eta}_{i}\mathbf{u}_{j}^{n+1}\cdot\tau_{i},\mathbf{v}\cdot\tau_{i}\rangle\nonumber\\
	&&\hspace{8mm}=(\mathbf{f}_j,\mathbf{v})_{\Omega_f}-\sum_{i=1}^{d-1}\langle(\eta_{i j}-\overline{\eta}_{i})\mathbf{u}_{j}^n\cdot\tau_{i},\mathbf{v}\cdot\tau_{i}\rangle+\langle\delta_{fj}^n,\mathbf{v}\cdot\mathbf{n}_{f}\rangle, \\
	&&b_f(\mathbf{u}_j^{n+1},q)=0,\label{decoupled-2}\\
	&&\gamma_{p}(\overline{\mathbb{K}}\triangledown\phi^{n+1}_{j},\triangledown\psi)_{\Omega_{p}}+g\langle\phi^{n+1}_{j},\psi\rangle
	=\langle\delta^n_{pj},\psi\rangle-\gamma_{p}((\mathbb{K}_{j}-\overline{\mathbb{K}})\triangledown\phi^n_{j},\triangledown\psi)_{\Omega_{p}}.\label{decoupled-3}
\end{eqnarray}

3. Update $\delta_{fj}^{n+1}$ and $\delta_{pj}^{n+1}$
\begin{eqnarray*}
&&\delta_{fj}^{n+1} =a\delta_{pj}^{n}+bg\phi_{j}^{n+1},\label{decoupled-comp1}\\
&&\delta_{pj}^{n+1} =c\delta_{fj}^{n}+d\mathbf{u}_{j}^{n+1}\cdot\mathbf{n}_f,\label{decoupled-comp3}
\end{eqnarray*}
where the coefficients
\begin{eqnarray*}
	a=\frac{\gamma_{f}}{\gamma_{p}},\hspace{9mm} b=-1-\frac{\gamma_{f}}{\gamma_{p}},\hspace{9mm}c=-1,\hspace{9mm} d=\gamma_{f}+\gamma_{p}.
\end{eqnarray*}

\section{ The Convergence of Parallel Ensemble Domain Decomposition Method} In this section, we focus on the case where the Robin parameters satisfy $\gamma_{f}\leq\gamma_{p}$. More importantly, we also prove that the proposed algorithm has the mesh-independent convergence rate in the case $\gamma_{f}<\gamma_{p}$. First, we would like to give an important Lemma \ref{a1b1c1<a2b2c2} which will be helpful in subsequent proofs.
\begin{lemma}\label{a1b1c1<a2b2c2}
	Assume that the parameters $a_1, a_2, b_1, b_2, c_1, c_2, d_1, d_2$ are positive constants, and $a_2<a_1, b_2<b_1, c_2<c_1, d_2<d_1$. For the iterative sequences $A^n,B^n,C^n,D^n,~n=1,2,\cdots$, if $a_1A^n+b_1B^n+c_1C^n+d_1D^n\leq a_2A^{n-1}+b_2B^{n-1}+c_2C^{n-1}+d_2D^{n-1}$,
	we have
	\begin{eqnarray*}
		a_1A^n+b_1B^n+c_1C^n+d_1D^n\leq \max\Big\{\frac{a2}{a1}, \frac{b2}{b1}, \frac{c2}{c1}, \frac{d2}{d1}\Big\}^{n-1} (a_2A^0+b_2B^0+c_2C^0+d_2D^0).
	\end{eqnarray*}
\end{lemma}

\begin{proof}
	For details, please refer to \cite{SunF}, which only adds one more item but proves similar.
\end{proof}

Moreover, before introducing the convergence theorem of the proposed algorithm, some quantities need to be defined
\begin{eqnarray*}
	&& \eta^{\prime\mathrm{max}}_{i j}=\max_{\mathbf{x}\in\Gamma} |{\eta_{i j}}(\mathbf{x})-\overline{\eta}_i(\mathbf{x})|,   \hspace{7mm} 
	\eta^{\prime\mathrm{max}}_{i}=\max_{j}\eta_{i j}^{\prime\mathrm{max}}, \\ 
	&& k_{j\min}=\min_{\mathbf{x}\in\Omega_p}{k_{j\min}(\mathbf{x})},\hspace{15mm} 
	k_{\min}=\min_{j}k_{j\min},\hspace{15mm}
	\overline{k}_{\min}=\min_{\mathbf{x}\in\Omega_p}\overline{k}_{\min}(\mathbf{x}),\\ 
	&& \rho_{j\max}^{\prime}=\max_{\mathbf{x}\in\Omega_p}\rho_{j}^{\prime}(\mathbf{x}),\hspace{20mm}  \rho^{\prime}_{\max}=\max_{j}\rho^{\prime}_{j\max}.
\end{eqnarray*}

\begin{theorem}
	Assume the solution of Parallel Ensemble DDM Algorithm (\ref{decoupled-1})-(\ref{decoupled-3}) is $(\mathbf{u}_{j}^n,p_{j}^n;\phi_{j}^n)$. Let $(\mathbf{u}_{j},p_{j};\phi_{j})$ denote the solution of the DDM weak formulation (\ref{decoupled-111})-(\ref{decoupled-222}),(\ref{decoupled-333}). Then if $\gamma_{f}\leq\gamma_{p}$ and $\delta_{fj}^{n}$, $\delta_{pj}^{n}$ satisfy the compatibility conditions (\ref{comp1})-(\ref{comp2}), assuming that $\overline{\eta}_{i}>\eta_{i}^{\prime \max}$ and $\overline{k}_{\min}>\rho_{\max}^{\prime}$, we can derive that $(\mathbf{u}_{j}^n,p_{j}^n;\phi_{j}^n)$ will converge to $(\mathbf{u}_{j},p_{j};\phi_{j})$.
\end{theorem}

\begin{proof}
Define the following notations for the error functions
\begin{eqnarray*}
&&\mathbf{e}^n_{uj}=\mathbf{u}_{j}-\mathbf{u}^n_{j},\hspace{9mm}{e}^n_{pj}=p_{j}-p^n_{j},\hspace{10mm} {e}^n_{\phi j}=\phi_j-\phi^n_j,\\
&& 
{\varepsilon}_{fj}^{n}=\delta_{fj}-{\delta}_{fj}^{n},\hspace{8mm}{\varepsilon}_{pj}^{n}=\delta_{pj}-{\delta}_{pj}^{n}.
\end{eqnarray*}
Then, for all $(\mathbf{v},q;\psi)\in (\mathbf{X}_{f},Q_f;X_p)$, by subtracting(\ref{decoupled-1})-(\ref{decoupled-3}) from
(\ref{decoupled-111})-(\ref{decoupled-222}),(\ref{decoupled-333}), we can have
\begin{eqnarray}\label{err-1} 
	&&a_{f}(\mathbf{e}_{uj}^{n+1},\mathbf{v})+b_f({e}_{pj}^{n+1},\mathbf{v})+\gamma_{f}\langle\mathbf{e}_{uj}^{n+1}\cdot\mathbf{n}_{f},\mathbf{v}\cdot\mathbf{n}_f\rangle+\sum_{i=1}^{d-1}\langle\overline{\eta}_{i}\mathbf{e}_{uj}^{n+1}\cdot\tau_{i},\mathbf{v}\cdot\tau_{i}\rangle\nonumber\\
    &&\hspace*{8mm}=-\sum_{i=1}^{d-1}\langle(\eta_{i j}-\overline{\eta}_{i})\mathbf{e}_{uj}^n\cdot\tau_{i},\mathbf{v}\cdot\tau_{i}\rangle+\langle{\varepsilon}_{fj}^n,\mathbf{v}\cdot\mathbf{n}_{f}\rangle,\\
    &&b_f(\mathbf{e}_{uj}^{n+1},q)=0,\label{err-2}\\
    &&\gamma_{p}(\overline{\mathbb{K}}\triangledown{e}^{n+1}_{\phi j},\triangledown\psi)_{\Omega_{p}}+g\langle {e}^{n+1}_{\phi j},\psi\rangle
    =\langle{\varepsilon}_{pj}^n,\psi\rangle-\gamma_{p}((\mathbb{K}_{j}-\overline{\mathbb{K}})\triangledown{e}^n_{\phi j},\triangledown\psi)_{\Omega_{p}}.\label{err-3}
\end{eqnarray}
Along the interface $\Gamma$, the error functions can be updated
\begin{eqnarray}
	&&\varepsilon_{fj}^{n+1} =a\varepsilon_{pj}^{n}+bg{e}_{\phi j}^{n+1},\label{errf}\\
	&&\varepsilon_{pj}^{n+1} =c\varepsilon_{fj}^{n}+d\mathbf{e}_{uj}^{n+1}\cdot\mathbf{n}_f.\label{errp}
\end{eqnarray}
Equation (\ref{errp}) can lead to
\begin{eqnarray}\label{etaD2}
||\varepsilon_{pj}^{n+1}||^2_{\Gamma} =c^2||\varepsilon_{fj}^{n}||^2_{\Gamma}+d^2||\mathbf{e}_{uj}^{n+1}\cdot\mathbf{n}_f||^2_{\Gamma}+2cd\langle\varepsilon_{fj}^{n},\mathbf{e}_{uj}^{n+1}\cdot\mathbf{n}_f\rangle.
\end{eqnarray}
Choosing $\mathbf{v}=\mathbf{e}_{uj}^{n+1}$ and $q=e_{pj}^{n+1}$ in (\ref{err-1})-(\ref{err-2}), we deduce
\begin{eqnarray}\label{aS+bS}
&&a_{f}(\mathbf{e}_{uj}^{n+1},\mathbf{e}_{uj}^{n+1})+\gamma_{f}\langle\mathbf{e}_{uj}^{n+1}\cdot\mathbf{n}_{f},\mathbf{e}_{uj}^{n+1}\cdot\mathbf{n}_f\rangle+\sum_{i=1}^{d-1}\langle\overline{\eta}_{i}\mathbf{e}_{uj}^{n+1}\cdot\tau_{i},\mathbf{e}_{uj}^{n+1}\cdot\tau_{i}\rangle\nonumber\\
&&\hspace{8mm}=-\sum_{i=1}^{d-1}\langle(\eta_{i j}-\overline{\eta}_{i})\mathbf{e}_{uj}^n\cdot\tau_{i},\mathbf{e}_{uj}^{n+1}\cdot\tau_{i}\rangle+\langle{\varepsilon}_{fj}^n,\mathbf{e}_{uj}^{n+1}\cdot\mathbf{n}_{f}\rangle.
\end{eqnarray}
Combining (\ref{etaD2}) and (\ref{aS+bS}), we can obtain
\begin{eqnarray}\label{errorD}
&&||\varepsilon_{pj}^{n+1}||^2_{\Gamma}=c^2||\varepsilon_{fj}^{n}||^2_{\Gamma}+(d^2+2cd\gamma_{f})||\mathbf{e}_{uj}^{n+1}\cdot\mathbf{n}_f||^2_{\Gamma}+2cda_{f}(\mathbf{e}_{uj}^{n+1},\mathbf{e}_{uj}^{n+1})\nonumber\\
&&\hspace*{18mm}+2cd\Big[\sum_{i=1}^{d-1}\langle\overline{\eta}_{i}\mathbf{e}_{uj}^{n+1}\cdot\tau_{i},\mathbf{e}_{uj}^{n+1}\cdot\tau_{i}\rangle+\sum_{i=1}^{d-1}\langle(\eta_{i j}-\overline{\eta}_{i})\mathbf{e}_{uj}^n\cdot\tau_{i},\mathbf{e}_{uj}^{n+1}\cdot\tau_{i}\rangle\Big].
\end{eqnarray}
Similarly,
\begin{eqnarray}\label{etaDD}
||\varepsilon_{fj}^{n+1}||^2_{\Gamma}
= a^2||\varepsilon_{pj}^{n}||^2_{\Gamma}+b^2g^2||{e}_{\phi j}^{n+1}||^2_{\Gamma}+2abg\langle\varepsilon_{pj}^{n},{e}_{\phi j}^{n+1}\rangle.
\end{eqnarray}
Set $\psi=g{e}_{\phi j}^{n+1}$ in (\ref{err-3}) to get
\begin{eqnarray}\label{ap+bp}
    g\langle{\varepsilon}_{pj}^n,{e}_{\phi j}^{n+1}\rangle=g\gamma_{p}(\overline{\mathbb{K}}\triangledown{e}^{n+1}_{\phi j},\triangledown {e}_{\phi j}^{n+1})_{\Omega_{p}}+g\gamma_{p}((\mathbb{K}_{j}-\overline{\mathbb{K}})\triangledown{e}^n_{\phi j},\triangledown{e}_{\phi j}^{n+1})_{\Omega_{p}}+g^2\langle {e}^{n+1}_{\phi j},{e}_{\phi j}^{n+1}\rangle.\hspace{7mm}
\end{eqnarray}
Combine (\ref{etaDD}) and (\ref{ap+bp}), we can address
\begin{eqnarray}\label{errorS}
	||\varepsilon_{fj}^{n+1}||^2_{\Gamma}&=&a^2||\varepsilon_{pj}^{n}||^2_{\Gamma}+(b^2+2ab)g^2||{e}_{\phi j}^{n+1}||^2_{\Gamma}\nonumber\\
	&&+2abg\gamma_{p}\Big[(\overline{\mathbb{K}}\triangledown{e}^{n+1}_{\phi j},\triangledown{e}_{\phi j}^{n+1})_{\Omega_{p}}+((\mathbb{K}_{j}-\overline{\mathbb{K}})\triangledown{e}^n_{\phi j},\triangledown{e}_{\phi j}^{n+1})_{\Omega_{p}}\Big],
\end{eqnarray}
since $a=\frac{\gamma_{f}}{\gamma_{p}}, b=-1-\frac{\gamma_{f}}{\gamma_{p}}, c=-1, d=\gamma_{f}+\gamma_{p}$, we summarize the following two important equations, which are the key point of convergence analysis.
\begin{eqnarray}\label{ffpp1}
	||\varepsilon_{pj}^{n+1}||^2_{\Gamma}&=&||\varepsilon_{fj}^{n}||^2_{\Gamma}+(\gamma_{p}^2-\gamma_{f}^2)||\mathbf{e}_{uj}^{n+1}\cdot\mathbf{n}_f||^2_{\Gamma}-2(\gamma_{f}+\gamma_{p})\sum_{i=1}^{d-1}\langle\overline{\eta}_{i}\mathbf{e}_{uj}^{n+1}\cdot\tau_{i},\mathbf{e}_{uj}^{n+1}\cdot\tau_{i}\rangle\nonumber\\
	&&-2(\gamma_{f}+\gamma_{p})a_{f}(\mathbf{e}_{uj}^{n+1},\mathbf{e}_{uj}^{n+1})-2(\gamma_{f}+\gamma_{p})\sum_{i=1}^{d-1}\langle(\eta_{i j}-\overline{\eta}_{i})\mathbf{e}_{uj}^n\cdot\tau_{i},\mathbf{e}_{uj}^{n+1}\cdot\tau_{i}\rangle,\\
	||\varepsilon_{fj}^{n+1}||^2_{\Gamma}&=&(\frac{\gamma_{f}}{\gamma_{p}})^2||\varepsilon_{pj}^{n}||^2_{\Gamma}+g^2(1-(\frac{\gamma_{f}}{\gamma_{p}})^2)||{e}_{\phi j}^{n+1}||^2_{\Gamma}\nonumber\\
	&&-2g\gamma_{f}(1+\frac{\gamma_{f}}{\gamma_{p}})\Big[(\overline{\mathbb{K}}\triangledown{e}^{n+1}_{\phi j},\triangledown{e}_{\phi j}^{n+1})_{\Omega_{p}}+((\mathbb{K}_{j}-\overline{\mathbb{K}})\triangledown{e}^n_{\phi j},\triangledown{e}_{\phi j}^{n+1})_{\Omega_{p}}\Big].\label{ffpp2}
\end{eqnarray}
Due to the convergence analysis being different between $\gamma_{f}=\gamma_{p}$ and $\gamma_{f}<\gamma_{p}$, we will discuss the two cases separately. \newline
\textbf{Case 1:$\gamma_{f}=\gamma_{p}=\gamma$} \newline 
For this case, we can rewrite (\ref{ffpp1})-(\ref{ffpp2}) as follows
\begin{eqnarray}\label{ffpp3}
	&&||\varepsilon_{pj}^{n+1}||^2_{\Gamma}=||\varepsilon_{fj}^{n}||^2_{\Gamma}-4\gamma a_{f}(\mathbf{e}_{uj}^{n+1},\mathbf{e}_{uj}^{n+1})-4\gamma\sum_{i=1}^{d-1}\langle\overline{\eta}_{i}\mathbf{e}_{uj}^{n+1}\cdot\tau_{i},\mathbf{e}_{uj}^{n+1}\cdot\tau_{i}\rangle\nonumber\\
	&&\hspace*{20mm}-4\gamma\sum_{i=1}^{d-1}\langle(\eta_{i j}-\overline{\eta}_{i})\mathbf{e}_{uj}^n\cdot\tau_{i},\mathbf{e}_{uj}^{n+1}\cdot\tau_{i}\rangle,\\
	&&||\varepsilon_{fj}^{n+1}||^2_{\Gamma}
	=||\varepsilon_{pj}^{n}||^2_{\Gamma}-4g\gamma(\overline{\mathbb{K}}\triangledown{e}^{n+1}_{\phi j},\triangledown{e}_{\phi j}^{n+1})_{\Omega_{p}}-4g\gamma((\mathbb{K}_{j}-\overline{\mathbb{K}})\triangledown{e}^n_{\phi j},\triangledown{e}_{\phi j}^{n+1})_{\Omega_{p}}.\label{ffpp4}
\end{eqnarray}
Adding (\ref{ffpp3}) and (\ref{ffpp4}), and summing over $n$ from $n=0$ to $N$, we will get
\begin{eqnarray}\label{ffpp2+3}
	&&||\varepsilon_{pj}^{N+1}||^2_{\Gamma}+||\varepsilon_{fj}^{N+1}||^2_{\Gamma}
	=||\varepsilon_{fj}^{0}||^2_{\Gamma}+||\varepsilon_{pj}^{0}||^2_{\Gamma}-4\gamma \sum_{n=0}^{N}a_{f}(\mathbf{e}_{uj}^{n+1},\mathbf{e}_{uj}^{n+1})\nonumber\\
	&&\hspace{35mm}-4\gamma\sum_{n=0}^{N}\big[\sum_{i=1}^{d-1}\langle\overline{\eta}_{i}\mathbf{e}_{uj}^{n+1}\cdot\tau_{i},\mathbf{e}_{uj}^{n+1}\cdot\tau_{i}\rangle+\sum_{i=1}^{d-1}\langle(\eta_{i j}-\overline{\eta}_{i})\mathbf{e}_{uj}^n\cdot\tau_{i},\mathbf{e}_{uj}^{n+1}\cdot\tau_{i}\rangle\nonumber\\
	&&\hspace{35mm}+g(\overline{\mathbb{K}}\triangledown{e}^{n+1}_{\phi j},\triangledown{e}_{\phi j}^{n+1})_{\Omega_{p}}+g((\mathbb{K}_{j}-\overline{\mathbb{K}})\triangledown{e}^n_{\phi j},\triangledown{e}_{\phi j}^{n+1})_{\Omega_{p}}\big].
\end{eqnarray}
By the Cauchy-Schwarz inequality, we have
\begin{eqnarray}\label{cauchy_schwarzf}
	&&-\sum_{i=1}^{d-1}\langle(\eta_{i j}-\overline{\eta}_{i})\mathbf{e}_{uj}^n\cdot\tau_{i},\mathbf{e}_{uj}^{n+1}\cdot\tau_{i}\rangle
	\leq\sum_{i=1}^{d-1}\int_{\varGamma}|\eta_{i j}-\overline{\eta}_{i}||(\mathbf{e}_{uj}^n\cdot\tau_{i})(\mathbf{e}_{uj}^{n+1}\cdot\tau_{i})|ds\nonumber\\
	&&\leq\sum_{i=1}^{d-1}\int_{\varGamma}|\eta_{i j}^{\prime \max}||(\mathbf{e}_{uj}^n\cdot\tau_{i})(\mathbf{e}_{uj}^{n+1}\cdot\tau_{i})|ds\leq\sum_{i=1}^{d-1}\Big[\frac{\eta_{i}^{\prime \max}}{2}\int_{\varGamma}(\mathbf{e}_{uj}^n\cdot\tau_{i})^2ds+\frac{\eta_{i}^{\prime \max}}{2}\int_{\varGamma}(\mathbf{e}_{uj}^{n+1}\cdot\tau_{i})^2ds\Big]\nonumber\\
	&&=\sum_{i=1}^{d-1}\Big[\frac{\eta_{i}^{\prime \max}}{2}||\mathbf{e}_{uj}^n\cdot\tau_{i}||_{\varGamma}^2+\frac{\eta_{i}^{\prime \max}}{2}||\mathbf{e}_{uj}^{n+1}\cdot\tau_{i}||_{\varGamma}^2\Big].
\end{eqnarray}
Meanwhile, we can further utilize the $2$-norm of vectors and matrices to analyze
\begin{eqnarray}\label{cauchy_schwarzp}
	&&-g((\mathbb{K}_{j}-\overline{\mathbb{K}})\triangledown e^n_{\phi j},\triangledown e_{\phi j}^{n+1})_{\Omega_{p}}\leq g\int_{\Omega_{p}}||\mathbb{K}_{j}-\overline{\mathbb{K}}||_2||\triangledown e^n_{\phi j}||_2||\triangledown e^{n+1}_{\phi j}||_2dx\nonumber\\	
	&&\leq g\int_{\Omega_{p}}\rho_j^{\prime}(x)||\triangledown e^n_{\phi j}||_2||\triangledown e^{n+1}_{\phi j}||_2dx\leq g\rho_{j \max}^{\prime}\int_{\Omega_{p}}||\triangledown e^n_{\phi j}||_2||\triangledown e^{n+1}_{\phi j}||_2dx\nonumber\\
	&&\leq g\rho_{j \max}^{\prime}||\triangledown e^n_{\phi j}||_{\Omega_{p}}||\triangledown e^{n+1}_{\phi j}||_{\Omega_{p}}\leq\frac{g}{2}\rho_{\max}^{\prime}||\triangledown e^n_{\phi j}||_{\Omega_{p}}^2+\frac{g}{2}\rho_{\max}^{\prime}||\triangledown e^{n+1}_{\phi j}||_{\Omega_{p}}^2.
\end{eqnarray}
Combining (\ref{cauchy_schwarzf})-(\ref{cauchy_schwarzp}) with (\ref{ffpp2+3}), the derivation is that
\begin{eqnarray}\label{ffpp5}
	&&||\varepsilon_{pj}^{N+1}||^2_{\Gamma}+||\varepsilon_{fj}^{N+1}||^2_{\Gamma}\nonumber\\
	&&\leq||\varepsilon_{fj}^{0}||^2_{\Gamma}+||\varepsilon_{pj}^{0}||^2_{\Gamma}+2\gamma\sum_{n=0}^{N}\sum_{i=1}^{d-1}(\eta_{i}^{\prime \max}||\mathbf{e}_{uj}^n\cdot\tau_{i}||_{\varGamma}^2-(2\overline{\eta}_{i}-\eta_{i}^{\prime \max})||\mathbf{e}_{uj}^{n+1}\cdot\tau_{i}||_{\varGamma}^2)\nonumber\\
	&&\hspace{3mm}-4\gamma \sum_{n=0}^{N}a_{f}(\mathbf{e}_{uj}^{n+1},\mathbf{e}_{uj}^{n+1})+2g\gamma\sum_{n=0}^{N}\big[\rho_{\max}^{\prime}||\triangledown{e}^n_{\phi j}||_{\Omega_{p}}^2-(2\overline{k}_{\min}-\rho_{\max}^{\prime})||\triangledown{e}^{n+1}_{\phi j}||_{\Omega_{p}}^2\big]\nonumber\\
	&&\leq||\varepsilon_{fj}^{0}||^2_{\Gamma}+||\varepsilon_{pj}^{0}||^2_{\Gamma}\nonumber+2\gamma\sum_{n=0}^{N+1}\sum_{i=1}^{d-1}\eta_{i}^{\prime \max}||\mathbf{e}_{uj}^n\cdot\tau_{i}||_{\varGamma}^2-2\gamma\sum_{n=0}^{N}\sum_{i=1}^{d-1}(2\overline{\eta}_{i}-\eta_{i}^{\prime \max})||\mathbf{e}_{uj}^{n+1}\cdot\tau_{i}||_{\varGamma}^2\nonumber\\
	&&\hspace{3mm}-4\gamma \sum_{n=0}^{N}a_{f}(\mathbf{e}_{uj}^{n+1},\mathbf{e}_{uj}^{n+1})+2g\gamma\big[\sum_{n=0}^{N+1}\rho_{\max}^{\prime}||\triangledown{e}^n_{\phi j}||_{\Omega_{p}}^2-\sum_{n=0}^{N}(2\overline{k}_{\min}-\rho_{\max}^{\prime})||\triangledown{e}^{n+1}_{\phi j}||_{\Omega_{p}}^2\big]\nonumber\\	&&\leq||\varepsilon_{fj}^{0}||^2_{\Gamma}+||\varepsilon_{pj}^{0}||^2_{\Gamma}+2\gamma\sum_{i=1}^{d-1}\eta_{i}^{\prime \max}||\mathbf{e}_{uj}^0\cdot\tau_{i}||_{\varGamma}^2+2g\gamma\rho_{\max}^{\prime}||\triangledown{e}^0_{\phi j}||_{\Omega_{p}}^2-4\gamma\sum_{n=0}^{N}a_{f}(\mathbf{e}_{uj}^{n+1},\mathbf{e}_{uj}^{n+1})\nonumber\\
	&&\hspace{3mm}-4\gamma\sum_{n=0}^{N}\sum_{i=1}^{d-1}(\overline{\eta}_{i}-\eta_{i}^{\prime \max})||\mathbf{e}_{uj}^{n+1}\cdot\tau_{i}||_{\varGamma}^2-4g\gamma\sum_{n=0}^{N}(\overline{k}_{\min}-\rho_{\max}^{\prime})||\triangledown{e}^{n+1}_{\phi j}||_{\Omega_{p}}^2.
\end{eqnarray}
Assume $\overline{\eta}_{i}-\eta_{i}^{\prime \max}>0$ and $\overline{k}_{\min}-\rho_{\max}^{\prime}>0$, then we have
\begin{eqnarray}\label{ffpp6}
	4\gamma \sum_{n=0}^{N}a_{f}(\mathbf{e}_{uj}^{n+1},\mathbf{e}_{uj}^{n+1})+2\gamma\sum_{n=0}^{N}\sum_{i=1}^{d-1}(\overline{\eta}_{i}-\eta_{i}^{\prime \max})||\mathbf{e}_{uj}^{n+1}\cdot\tau_{i}||_{\varGamma}^2+4g\gamma\sum_{n=0}^{N}(\overline{k}_{\min}-\rho_{\max}^{\prime})||\triangledown{e}^{n+1}_{\phi j}||_{\Omega_{p}}^2\nonumber\\
	\leq||\varepsilon_{fj}^{0}||^2_{\Gamma}+||\varepsilon_{pj}^{0}||^2_{\Gamma}+2\gamma\sum_{i=1}^{d-1}\eta_{i}^{\prime \max}||\mathbf{e}_{uj}^0\cdot\tau_{i}||_{\varGamma}^2+2g\gamma\rho_{\max}^{\prime}||\triangledown{e}^0_{\phi j}||_{\Omega_{p}}^2,
\end{eqnarray}
which means that $\mathbf{e}_{uj}^{n+1}, e^{n+1}_{\phi j}$ are convergent in the sense of $H^1(\Omega_f)^d$ and $H^1(\Omega_p)$ while $n\rightarrow \infty$. 

For the convergence of $\varepsilon^n_{pj}$, let $\psi\in H_{00}^\frac{1}{2}(\Gamma)$, $\theta_{\psi} \in X_p$ be a harmonic extension of $\psi$ to the porous media domain, which satisfies
\begin{eqnarray*}
	\begin{cases}
		-\Delta\theta_{\psi}=0 & \quad\text{in}~\Omega_{p}, \\
		\theta_{\psi}=\psi & \quad\text{on}~\Gamma, \\
		\theta_{\psi}=0 & \quad\text{on}~\partial\Omega_{p}/\Gamma,
	\end{cases}
\end{eqnarray*}
where $H_{00}^{\frac{1}{2}}(\Gamma)$ is the interpolation space
\begin{eqnarray*}
H_{00}^{\frac{1}{2}}(\Gamma)=[L^2(\Gamma),H_0^1(\Gamma)]_{\frac{1}{2}}.
\end{eqnarray*}
There exist the positive constants $C_{\theta1}, C_{\theta2}$, 
\begin{eqnarray*}
	||\theta_{\psi}||_{H^1(\Omega_{p})}\leq C_{\theta1}||\psi||_{H_{00}^{\frac{1}{2}}(\Gamma)}\leq C_{\theta2}||\psi||_{H^1(\Omega_{p})}.
\end{eqnarray*}
According to (\ref{err-3}), there exist positive constants $C_{\psi 1},C_{\psi 2},C_{\psi 3}$ such that
\begin{eqnarray}\label{convergence1}
	&&||\varepsilon^n_{pj}||_{H^{-\frac{1}{2}}(\Gamma)}=\sup\limits_{\forall\psi\neq0\in H^{\frac{1}{2}}_{00}(\Gamma)} \frac{\langle\varepsilon^n_{pj},\psi\rangle}{||\psi||_{H_{00}^{\frac{1}{2}}(\Gamma)}}\leq\sup\limits_{\forall\theta_{\psi}\neq0\in X_p}\frac{\langle\varepsilon^n_{pj},\theta_{\psi}\rangle}{||\theta_{\psi}||_{H^1(\Omega_p)}}\nonumber\\
	&&\leq\sup\limits_{\forall\theta_{\psi}\neq0\in X_p}\frac{\gamma(\overline{\mathbb{K}}\triangledown{e}^{n+1}_{\phi j},\triangledown\theta_{\psi})_{\Omega_{p}}}{||\theta_{\psi}||_{H^1(\Omega_p)}}+\sup\limits_{\forall\theta_{\psi}\neq0\in X_p}\frac{\gamma((\mathbb{K}_{j}-\overline{\mathbb{K}})\triangledown{e}^{n}_{\phi j},\triangledown\theta_{\psi})_{\Omega_{p}}}{||\theta_{\psi}||_{H^1(\Omega_p)}}+\sup\limits_{\forall\theta_{\psi}\neq0\in X_p}\frac{g\langle e^{n+1}_{\phi j},\theta_{\psi}\rangle}{||\theta_{\psi}||_{H^1(\Omega_p)}}\nonumber\\
	&&\leq(C_{\psi 1}\gamma\overline{k}_{\max}+C_{\psi 3}g)||\triangledown{e}^{n+1}_{\phi j}||_{H^1(\Omega_p)}+C_{\psi 2}\gamma\rho_{\max}^{\prime}||\triangledown{e}^{n}_{\phi j}||_{H^1(\Omega_p)}.
\end{eqnarray}
Besides, according to (\ref{errf}), we have
\begin{eqnarray}\label{convergence2}
	||\varepsilon^{n}_{fj}||_{H^{-\frac{1}{2}}(\Gamma)}
	&\leq&\frac{\gamma_{f}}{\gamma_{p}}||\varepsilon^{n-1}_{pj}||_{H^{-\frac{1}{2}}(\Gamma)}+(\gamma_{f}+\gamma_{p})||g{e}^{n}_{\phi j}||_{H^{-\frac{1}{2}}(\Gamma)}\nonumber\\
	&\leq&||\varepsilon^{n-1}_{pj}||_{H^{-\frac{1}{2}}(\Gamma)}+2\gamma||g{e}^{n}_{\phi j}||_{H^1(\Omega_p)}.
\end{eqnarray}
Thus we can get the convergence of $\varepsilon^n_{pj},\varepsilon^n_{fj}$ in $H^{-\frac{1}{2}}(\Gamma)$. 
Then, there exist positive constants $C_{v 1},C_{v 2},C_{v 3},C_{v 4},C_{v 5}$ such that 
	\begin{eqnarray}\label{pshoulian}
		&&||e^n_{pj}||_{\Omega_{f}}\leq\sup\limits_{\forall\mathbf{v}\neq0\in\mathbf{X}_{f}}\frac{(e^n_{pj},\triangledown\cdot\mathbf{v})}{||\mathbf{v}||_{H^1(\Omega_{f})}}\nonumber\\
		&&\leq\sup\limits_{\forall\mathbf{v}\neq0\in\mathbf{X}_{f}}\frac{a_{f}(\mathbf{e}_{uj}^{n},\mathbf{v})}{||\mathbf{v}||_{H^1(\Omega_{f})}}+\sup\limits_{\forall\mathbf{v}\neq0\in\mathbf{X}_{f}}\frac{\gamma\langle\mathbf{e}_{uj}^{n}\cdot\mathbf{n}_{f},\mathbf{v}\cdot\mathbf{n}_f\rangle}{||\mathbf{v}||_{H^1(\Omega_{f})}}+\sup\limits_{\forall\mathbf{v}\neq0\in\mathbf{X}_{f}}\frac{\sum_{i=1}^{d-1}\langle\overline{\eta}_{i}\mathbf{e}_{uj}^{n}\cdot\tau_{i},\mathbf{v}\cdot\tau_{i}\rangle}{||\mathbf{v}||_{H^1(\Omega_{f})}}\nonumber\\
		&&\hspace{8mm}+\sup\limits_{\forall\mathbf{v}\neq0\in\mathbf{X}_{f}}\frac{\sum_{i=1}^{d-1}\langle(\eta_{i j}-\overline{\eta}_{i})\mathbf{e}_{uj}^{n-1}\cdot\tau_{i},\mathbf{v}\cdot\tau_{i}\rangle}{||\mathbf{v}||_{H^1(\Omega_{f})}}-\sup\limits_{\forall\mathbf{v}\neq0\in\mathbf{X}_{f}}\frac{\langle{\varepsilon}_{fj}^n,\mathbf{v}\cdot\mathbf{n}_{f}\rangle}{||\mathbf{v}||_{H^1(\Omega_{f})}}\nonumber\\
		&&\leq[C_{v 1}\nu+C_{v 2}\gamma+C_{v 3}\overline{\eta}_i]||\mathbf{e}_{uj}^{n}||_{H^1(\Omega_{f})}+C_{v 4}\eta^{\prime\mathrm{max}}_{i}||\mathbf{e}_{uj}^{n-1}||_{H^1(\Omega_{f})}+C_{v 5}||\varepsilon_{fj}^n||_{H^{-\frac{1}{2}}(\Gamma)},
\end{eqnarray}
which implies the convergence of the kinematic pressure $e^n_{pj}$.\newline
\textbf{Case 2 : $\gamma_{f}< \gamma_{p} $}\newline
In this case, thanks to the trace inequality and Poincar\'e inequality, there exists a positive constant $C_p$, such that
\begin{eqnarray}\label{ineq01}
	||{e}^{n+1}_{\phi j}||_{\varGamma}^2\leq C_p||\triangledown{e}^{n+1}_{\phi j}||^2_{\Omega_{p}}.
\end{eqnarray}
By (\ref{ffpp2}), (\ref{cauchy_schwarzp}) and the Cauchy-Schwarz inequality, we have
\begin{eqnarray}\label{ineq03}
	&&g^2\big[1-(\frac{\gamma_{f}}{\gamma_{p}})^2\big]||{e}_{\phi j}^{n+1}||^2_{\Gamma}-2g\gamma_{f}(1+\frac{\gamma_{f}}{\gamma_{p}})[(\overline{\mathbb{K}}\triangledown{e}^{n+1}_{\phi j},\triangledown{e}_{\phi j}^{n+1})_{\Omega_{p}}
	+((\mathbb{K}_{j}-\overline{\mathbb{K}})\triangledown{e}^n_{\phi j},\triangledown{e}_{\phi j}^{n+1})_{\Omega_{p}}]\nonumber\\
	&&\leq g^2\big[1-(\frac{\gamma_{f}}{\gamma_{p}})^2\big]C_p||\triangledown{e}^{n+1}_{\phi j}||^2_{\Omega_{p}}+g\gamma_{f}(1+\frac{\gamma_{f}}{\gamma_{p}})\rho_{\max}^{\prime}(||\triangledown{e}^n_{\phi j}||_{\Omega_{p}}^2+||\triangledown{e}^{n+1}_{\phi j}||_{\Omega_{p}}^2)\nonumber\\
	&&\hspace{3mm}-2g\gamma_{f}\big[1+\frac{\gamma_{f}}{\gamma_{p}}\big](\overline{\mathbb{K}}\triangledown{e}^{n+1}_{\phi j},\triangledown {e}_{\phi j}^{n+1})_{\Omega_{p}}\nonumber\\
    &&\leq g(1+\frac{\gamma_{f}}{\gamma_{p}})[g(1-\frac{\gamma_{f}}{\gamma_{p}})C_p-\gamma_{f}(2\overline{k}_{\min}-\rho_{\max}^{\prime})]||\triangledown{e}^{n+1}_{\phi j}||^2_{\Omega_{p}}\nonumber\\
    &&\hspace{3mm}+g\gamma_{f}(1+\frac{\gamma_{f}}{\gamma_{p}})\rho_{\max}^{\prime}||\triangledown{e}^n_{\phi j}||_{\Omega_{p}}^2.
\end{eqnarray}
Integrating the above inequalities with (\ref{ffpp2}), we obtain
\begin{eqnarray}\label{case201}
	||\varepsilon_{fj}^{n+1}||^2_{\Gamma}
	&\leq&(\frac{\gamma_{f}}{\gamma_{p}})^2||\varepsilon_{pj}^{n}||^2_{\Gamma}+g(1+\frac{\gamma_{f}}{\gamma_{p}})[g(1-\frac{\gamma_{f}}{\gamma_{p}})C_p-\gamma_{f}(2\overline{k}_{\min}-\rho_{\max}^{\prime})]||\triangledown{e}^{n+1}_{\phi j}||^2_{\Omega_{p}}\nonumber\\
	&&+g\gamma_{f}(1+\frac{\gamma_{f}}{\gamma_{p}})\rho_{\max}^{\prime}||\triangledown{e}^n_{\phi j}||_{\Omega_{p}}^2.
\end{eqnarray}
Similarly,
 there exists a positive constant $C_f$ such that
\begin{eqnarray}\label{ineq04}
	||\mathbf{e}_{uj}^{n+1}\cdot\mathbf{n}_f||^2_{\Gamma}\leq C_f||\mathbb{D}(\mathbf{e}_{uj}^{n+1})||_{\Omega_{f}}^2.
\end{eqnarray}
By utilizing the Cauchy-Schwarz inequality and (\ref{cauchy_schwarzf}),
 we also address
\begin{eqnarray}\label{ineq06}
	&&(\gamma_{p}^2-\gamma_{f}^2)||\mathbf{e}_{uj}^{n+1}\cdot\mathbf{n}_f||^2_{\Gamma}
	-2(\gamma_{f}+\gamma_{p})a_{f}(\mathbf{e}_{uj}^{n+1},\mathbf{e}_{uj}^{n+1})\nonumber-2(\gamma_{f}+\gamma_{p})\sum_{i=1}^{d-1}\overline{\eta}_{i}||\mathbf{e}_{uj}^{n+1}\cdot\tau_{i}||_{\varGamma}^2\nonumber\\
	&&\hspace{3mm}-2(\gamma_{f}+\gamma_{p})\sum_{i=1}^{d-1}\langle(\eta_{i j}-\overline{\eta}_{i})\mathbf{e}_{uj}^n\cdot\tau_{i},\mathbf{e}_{uj}^{n+1}\cdot\tau_{i}\rangle\nonumber\\
	&&\leq(\gamma_{p}^2-\gamma_{f}^2)C_f||\mathbb{D}(\mathbf{e}_{uj}^{n+1})||_{\Omega_{f}}^2-2(\gamma_{f}+\gamma_{p})a_{f}(\mathbf{e}_{uj}^{n+1},\mathbf{e}_{uj}^{n+1})-2(\gamma_{f}+\gamma_{p})\sum_{i=1}^{d-1}\overline{\eta}_{i}||\mathbf{e}_{uj}^{n+1}\cdot\tau_{i}||_{\varGamma}^2\nonumber\\
	&&\hspace{3mm}+(\gamma_{f}+\gamma_{p})\sum_{i=1}^{d-1}\big[\eta_{i}^{\prime \max}||\mathbf{e}_{uj}^n\cdot\tau_{i}||_{\varGamma}^2+\eta_{i}^{\prime \max}||\mathbf{e}_{uj}^{n+1}\cdot\tau_{i}||_{\varGamma}^2\big]\nonumber\nonumber\\
	&&\leq\big[(\gamma_{p}^2-\gamma_{f}^2)C_f-4\nu(\gamma_{f}+\gamma_{p})\big]||\mathbb{D}(\mathbf{e}_{uj}^{n+1})||_{\Omega_{f}}^2+(\gamma_{f}+\gamma_{p})\sum_{i=1}^{d-1}\eta_{i}^{\prime \max}||\mathbf{e}_{uj}^n\cdot\tau_{i}||_{\varGamma}^2\nonumber\\
	&&\hspace{3mm}-(\gamma_{f}+\gamma_{p})\sum_{i=1}^{d-1}(2\overline{\eta}_{i}-\eta_{i}^{\prime \max})||\mathbf{e}_{uj}^{n+1}\cdot\tau_{i}||_{\varGamma}^2.
\end{eqnarray}
Together with (\ref{ffpp2}), we can get
\begin{eqnarray}\label{case202}
	&&||\varepsilon_{pj}^{n+1}||^2_{\Gamma}
	\leq||\varepsilon_{fj}^{n}||^2_{\Gamma}+\big[(\gamma_{p}^2-\gamma_{f}^2)C_f-4\nu(\gamma_{f}+\gamma_{p})\big]||\mathbb{D}(\mathbf{e}_{uj}^{n+1})||_{\Omega_{f}}^2\nonumber\\
	&&\hspace{18mm}+(\gamma_{f}+\gamma_{p})\sum_{i=1}^{d-1}\big[\eta_{i}^{\prime \max}||\mathbf{e}_{uj}^n\cdot\tau_{i}||_{\varGamma}^2-(2\overline{\eta}_{i}-\eta_{i}^{\prime \max})||\mathbf{e}_{uj}^{n+1}\cdot\tau_{i}||_{\varGamma}^2\big].
\end{eqnarray}
Let
\begin{eqnarray}\label{casecondition2}
	\gamma_{p}-\gamma_{f}<\frac{4\nu}{C_f},
\end{eqnarray}
which leads to $\big[(\gamma_{p}^2-\gamma_{f}^2)C_f-4\nu(\gamma_{f}+\gamma_{p})\big]||\mathbb{D}(\mathbf{e}_{uj}^{n+1})||_{\Omega_{f}}^2<0$.

Multiply $\frac{1}{2}\big[1+(\frac{\gamma_{f}}{\gamma_{p}})^2\big]$ with (\ref{case202}) and combine with (\ref{case201}), we can get
\begin{eqnarray}\label{case203}
	&&\frac{1}{2}\big[1+(\frac{\gamma_{f}}{\gamma_{p}})^2\big]||\varepsilon_{pj}^{n+1}||^2_{\Gamma}+||\varepsilon_{fj}^{n+1}||^2_{\Gamma}+\frac{1}{2}\big[1+(\frac{\gamma_{f}}{\gamma_{p}})^2\big](\gamma_{f}+\gamma_{p})\sum_{i=1}^{d-1}(2\overline{\eta}_{i}-\eta_{i}^{\prime \max})||\mathbf{e}_{uj}^{n+1}\cdot\tau_{i}||_{\varGamma}^2\nonumber\\
	&&\hspace{8mm}+g(1+\frac{\gamma_{f}}{\gamma_{p}})\big[\gamma_{f}(2\overline{k}_{\min}-\rho_{\max}^{\prime})-g(1-\frac{\gamma_{f}}{\gamma_{p}})C_p\big]||\triangledown{e}^{n+1}_{\phi j}||^2_{\Omega_{p}}\nonumber\\
	&&\leq\frac{1}{2}\big[1+(\frac{\gamma_{f}}{\gamma_{p}})^2\big]||\varepsilon_{fj}^{n}||^2_{\Gamma}+\frac{1}{2}\big[1+(\frac{\gamma_{f}}{\gamma_{p}})^2\big](\gamma_{f}+\gamma_{p})\sum_{i=1}^{d-1}\eta_{i}^{\prime \max}||\mathbf{e}_{uj}^n\cdot\tau_{i}||_{\Gamma}^2+(\frac{\gamma_{f}}{\gamma_{p}})^2||\varepsilon_{pj}^{n}||^2_{\Gamma}\nonumber\\
	&&\hspace{8mm}+g\gamma_{f}(1+\frac{\gamma_{f}}{\gamma_{p}})\rho_{\max}^{\prime}||\triangledown{e}^n_{\phi j}||_{\Omega_{p}}^2.
\end{eqnarray}
We assume that
\begin{eqnarray}\label{casecondition1}
	\frac{1}{\gamma_{f}}-\frac{1}{\gamma_{p}}<\frac{2(\overline{k}_{\min}-\rho_{\max}^{\prime})}{gC_p},
\end{eqnarray} 
then, it's easy to check 
\begin{eqnarray*}
	&&a_2=(\frac{\gamma_{f}}{\gamma_{p}})^2<a_1=\frac{1}{2}\big[1+(\frac{\gamma_{f}}{\gamma_{p}})^2\big], \\
	&&b_2=\frac{1}{2}\big[1+(\frac{\gamma_{f}}{\gamma_{p}})^2\big]<b_1=1, \\
	&&c_2=\frac{1}{2}\big[1+(\frac{\gamma_{f}}{\gamma_{p}})^2\big](\gamma_{f}+\gamma_{p})\sum_{i=1}^{d-1}\eta_{i}^{\prime \max}<c_1=\frac{1}{2}[1+(\frac{\gamma_{f}}{\gamma_{p}})^2](\gamma_{f}+\gamma_{p})\sum_{i=1}^{d-1}(2\overline{\eta}_{i}-\eta_{i}^{\prime \max}),\\
	&&d_2=g\gamma_{f}(1+\frac{\gamma_{f}}{\gamma_{p}})\rho_{\max}^{\prime}<d_1=g(1+\frac{\gamma_{f}}{\gamma_{p}})\big[\gamma_{f}(2\overline{k}_{\min}-\rho_{\max}^{\prime})-g(1-\frac{\gamma_{f}}{\gamma_{p}})C_p\big].
\end{eqnarray*}
Hence, according to Lemma \ref{a1b1c1<a2b2c2}, we arrive at
\begin{eqnarray}
		&&\frac{1}{2}(1+(\frac{\gamma_{f}}{\gamma_{p}})^2)||\varepsilon_{pj}^{N}||^2_{\Gamma}+||\varepsilon_{fj}^{N}||^2_{\Gamma}+\frac{1}{2}\big[1+(\frac{\gamma_{f}}{\gamma_{p}})^2\big](\gamma_{f}+\gamma_{p})\sum_{i=1}^{d-1}(2\overline{\eta}_{i}-\eta_{i}^{\prime \max})||\mathbf{e}_{uj}^{N}\cdot\tau_{i}||_{\Gamma}^2\nonumber\\
&&\hspace{8mm}+g\big[1+(\frac{\gamma_{f}}{\gamma_{p}})\big]\big[\gamma_{f}(2\overline{k}_{\min}-\rho_{\max}^{\prime})-gC_p\big(1-(\frac{\gamma_{f}}{\gamma_{p}})\big)\big]||\triangledown {e}_{\phi j}^{N}||^2_{\Omega_{p}}\nonumber\\
&&\leq\mathcal{E}^{N-1}\bigg\{(\frac{\gamma_{f}}{\gamma_{p}})^2||\varepsilon_{pj}^{0}||^2_{\Gamma}+\frac{1}{2}\big[1+(\frac{\gamma_{f}}{\gamma_{p}})^2\big]||\varepsilon_{fj}^{0}||^2_{\Gamma}+g\gamma_{f}(1+\frac{\gamma_{f}}{\gamma_{p}})\rho_{\max}^{\prime}||\triangledown{e}^0_{\phi j}||_{\Omega_{p}}^2\nonumber\\
&&\hspace{8mm}+\frac{1}{2}\big[1+(\frac{\gamma_{f}}{\gamma_{p}})^2\big](\gamma_{f}+\gamma_{p})\sum_{i=1}^{d-1}\eta_{i}^{\prime \max}||\mathbf{e}_{uj}^0\cdot\tau_{i}||_{\Gamma}^2\bigg\},
\end{eqnarray}
where
\begin{eqnarray}\label{E}
	\mathcal{E}=\max\Big\{\frac{2(\frac{\gamma_{f}}{\gamma_{p}})^2}{1+(\frac{\gamma_{f}}{\gamma_{p}})^2}, \frac{1}{2}\big[1+(\frac{\gamma_{f}}{\gamma_{p}})^2\big], \sum_{i=1}^{d-1} \frac{\eta_{i}^{\prime \max}}{2\overline{\eta}_{i}-\eta_{i}^{\prime \max}},\frac{\rho_{\max}^{\prime}}{(2\overline{k}_{\min}-\rho_{\max}^{\prime})-gC_p(\frac{1}{\gamma_{f}}-\frac{1}{\gamma_{p}})}\Big\},\hspace{7mm}
\end{eqnarray}
which illustrates that $||\varepsilon_{pj}^{n}||^2_{\Gamma}$, $||\varepsilon_{fj}^{n}||^2_{\Gamma}$, $||\mathbf{e}_{uj}^{n}\cdot\tau_{i}||^2_{\Gamma}$ and $||\triangledown {e}_{\phi j}^{n}||^2_{\Omega_{p}}$ are mesh-indenpedent geometrically convergent. Then, $||\varepsilon_{fj}^{n}||^2_{\Gamma}$ and $||\mathbf{e}_{uj}^{n}\cdot\tau_{i}||^2_{\Gamma}$ together with the error equation (\ref{aS+bS}) further implies the geometric convergence of the velocity $||\triangledown\mathbf{e}_{uj}^n||_{\Omega_{f}}^2$. Besides, the convergence of the pressure $||e_{pj}^n||^2_{\Omega_{f}}$ can be derived by the similar way as (\ref{pshoulian}).
\end{proof}

Furthermore, we can derive the following mesh-independent geometric convergence result for the proposed algorithm while the Robin parameters satisfy $\gamma_{f}<\gamma_{p}$.
\begin{theorem}\label{case2<thm}
	For the case $\gamma_{f}<\gamma_{p}$, assume that
	\begin{eqnarray}
			\hspace{5mm}0<\gamma_{p}-\gamma_{f}<\frac{4\nu}{C_f},\hspace{9mm} 	\frac{1}{\gamma_{f}}-\frac{1}{\gamma_{p}}<\frac{(2\overline{k}_{\min}-2\rho_{\max}^{\prime})}{gC_p},\hspace{9mm} \overline{\eta}_{i}>\eta_{i}^{\prime \max},\hspace{9mm}\overline{k}_{\min}>\rho_{\max}^{\prime},\nonumber
	\end{eqnarray}
the parallel Ensemble DDM algorithm has the following mesh-independent geometric convergence
	\begin{eqnarray}\label{case2last}
		&&||\varepsilon_{pj}^{N}||^2_{\Gamma}
		+||e_{pj}^N||^2_{\Omega_{f}}+||\varepsilon_{fj}^{N}||^2_{\Gamma}+||\mathbf{e}_{uj}^N\cdot\tau_{i}||_{\Gamma}^2+||\triangledown {e}_{\phi j}^{N}||^2_{\Omega_{p}}+||\triangledown\mathbf{e}_{uj}^N||_{\Omega_{f}}^2\nonumber\\
		&&\leq C\mathcal{E}^{N-1}\bigg\{||\varepsilon_{pj}^{0}||^2_{\Gamma}+||\varepsilon_{fj}^{0}||^2_{\Gamma}+||\triangledown{e}^0_{\phi j}||_{\Omega_{p}}^2+||\mathbf{e}_{uj}^0\cdot\tau_{i}||_{\Gamma}^2\bigg\},
    \end{eqnarray}
where $\mathcal{E}$ is defined by (\ref{E}).
\end{theorem}


\begin{theorem}\label{remark}
	 The optimal Robin parameters $\gamma_{f}^{\ast}$ and $\gamma_{p}^{\ast}$ are
\begin{eqnarray}\label{opm}
	\gamma_{f}^{\ast}&=&\frac{1-2\mu_f|\overline{\mathbb{K}}|s_{\min}s_{ \max}}{|\overline{\mathbb{K}}|(s_{\min}+s_{\max})}+\sqrt{(\frac{1-2\mu_f|\overline{\mathbb{K}}|s_{\min}s_{\max}}{|\overline{\mathbb{K}|}(s_{\min}+s_{\max})})^2+\frac{2\mu_f}{|\overline{\mathbb{K}}|}}, \nonumber\\
	\gamma_{p}^{\ast}&=&-\frac{1-2\mu_f|\overline{\mathbb{K}}|s_{\min}s_{\max}}{|\overline{\mathbb{K}}|(s_{\min}+s_{\max})}+\sqrt{(\frac{1-2\mu_f|\overline{\mathbb{K}}|s_{\min}s_{\max}}{|\overline{\mathbb{K}}|(s_{\min}+s_{\max})})^2+\frac{2\mu_f}{|\overline{\mathbb{K}}|}},
\end{eqnarray}	
where $|\overline{\mathbb{K}}|$ means the determinant of the average of the hydraulic conductivity, and $ s_{\min}$,  $s_{\max}$ are the minimum and maximum values of the frequency variable $s$ of the Fourier transform. In general,  $s_{\min}=\frac{\pi}{L^{\prime}}$( $L^{\prime}$ is the length of the interface $\Gamma$) and $s_{\max}=\frac{\pi}{h}$( $h$ is the mesh size).
\end{theorem}
\begin{proof}
We refer to articles \cite{osm01} to select the optimal Robin parameters. Similarly, for $\omega_s(x,y)\in L^2(\mathbb{R}^2)$, we have the Fourier transform
\begin{eqnarray*}
\mathcal{F}:\omega_s(x,y)\mapsto\widehat{\omega}_s(x,s)=\int_{\mathbb{R}}e^{-isy}\omega_s(x,y)dy,
\end{eqnarray*}
where $s$ is the frequency variable. Inspired by \cite{osm01}, the relevant frequencies are assumed in a range of $ 0<s_{\min}\leq s\leq s_{\max}$.
 It's easy to find that (3.5) in the paper \cite{osm01} is in the same form as the problem we discussed, except for the third equation. So in this paper according to the characteristic of the problem, we need to change the third one as
\begin{eqnarray*}
-\nu(\partial_x\mathbf{e}_{u_2j}^m+\partial_y\mathbf{e}_{u_1j}^m)=\overline{\eta}\mathbf{e}_{u_2j}^m-(\overline{\eta}-\eta_j)\mathbf{e}_{u_2j}^{m-1}.
\end{eqnarray*}
where $\mathbf{e}_{u_1j}$ and $\mathbf{e}_{u_2j}$ represent the errors of the velocity components for Stokes part. This formula only affects $\mathbf{u}_2$ and has no influence on selecting the optimal Robin parameters.
Besides, assuming $\mathbb{K}_j=k_j\mathbb{I}$ and  $\overline{\mathbb{K}}=\overline{k}\mathbb{I}$, where $k_j$ and $\overline{k}$ are positive constants, the first expression of formula (3.6) in \cite{osm01} needs to be changed as
\begin{eqnarray}\label{rewrite}
-\Big(\partial_x(\overline{k}\partial_x)+\partial_y(\overline{k}\partial_y)\Big)e_{\phi j}^m=\Big(\partial_x((\overline{k}-k_j)\partial_x) +\partial_y((\overline{k}-k_j)\partial_y)\Big)e_{\phi j}^{m-1}.
\end{eqnarray}
Using the Fourier transform, the equation (\ref{rewrite}) yields the following transformation
\begin{eqnarray}\label{rewritepp}
	-\partial_x(\overline{k}\partial_x)\hat{e}_{\phi j}^m+s^2\overline{k}\hat{e}_{\phi j}^m&=&\partial_x((\overline{k}-k_j)\partial_x)\hat{e}_{\phi j}^{m-1}+s^2(\overline{k}-k_j)\hat{e}_{\phi j}^{m-1}\nonumber\\
	&=&\frac{\overline{k}-k_j}{\overline{k}}\Big(\partial_x(\overline{k}\partial_x)\hat{e}_{\phi j}^{m-1}+s^2\overline{k}\hat{e}_{\phi j}^{m-1}\Big)\nonumber\\
	&=&\frac{\overline{k}-k_j}{\overline{k}}\Big(\partial_x((\overline{k}-k_j)\partial_x)\hat{e}_{\phi j}^{m-2}+s^2(\overline{k}-k_j)\hat{e}_{\phi j}^{m-2}\Big)\nonumber\\
	&=&\cdots =\Big(\frac{\overline{k}-k_j}{\overline{k}}\Big)^{m-1}\Big(\partial_x((\overline{k}-k_j)\partial_x)\hat{e}_{\phi j}^{0}+s^2(\overline{k}-k_j)\hat{e}_{\phi j}^{0}\Big),\hspace{3mm}
\end{eqnarray}
where $\hat{e}_{\phi j}^0=\hat{\phi}_j-\hat{\phi}_j^0$. In the first step of ensemble DDM algorithm, we set $\hat{\phi}_j^0=0$. Since the sink/source term in (\ref{Darcy01}) is assumed zero, then the true solution $\hat{\phi}_j$ can be calculated directly by Fourier transform as \cite{osm01}. Then, we have
\begin{eqnarray}
	-\partial_x\overline{k}\partial_x\hat{e}_{\phi j}^m+s^2\overline{k}\hat{e}_{\phi j}^m=-\partial_x\overline{k}\partial_x\hat{\phi}_j+s^2\overline{k}\hat{\phi}_j=0.
\end{eqnarray}
 It can be inferred through analysis that the selection of the optimal Robin parameters in this paper is consistent with that in \cite{osm01}. Therefore, we can summarize the optimal Robin parameters $\gamma_{f}^{\ast}$ and $\gamma_{p}^{\ast}$ as (\ref{opm}).
\end{proof}
 
\section{ Finite Element Approximation}
In this section, we will propose the finite element discretization of the above parallel ensemble domain decomposition method for solving the random Stokes-Darcy coupled model. The advantage of considering finite element approximations is that we could obtain the explicit mesh-dependent convergence of such algorithm in the case  $\gamma_{f}=\gamma_{p}=\gamma$.

We prefer to consider a regular and quasi-uniform triangulation  $\mathcal{T}_h$ for the global domain $\overline{\Omega }_{f}\cup \overline{\Omega}_{p}$, with the property that $\mathcal{T}_{fh}$, $\mathcal{T}_{ph}$ are two triangulations of the subdomains $\Omega_{f},\Omega_{p}$. We should assume that $\mathcal{T}_{fh}$, $\mathcal{T}_{ph}$ are compatible on the interface $\Gamma$ and the mesh is quasi-uniform on $\Gamma$. We denote the introduced triangulation on $\Gamma$ as $\mathcal{B}_h$. We define the corresponding finite element space $\mathbf{X}_{fh} \subset \mathbf{X}_f, Q_{fh} \subset Q_f$ for $\Omega_{f}$, and $X_{ph} \subset X_p$ for $\Omega_{p}$.
\begin{eqnarray*}
	\mathbf{X}_{fh}&:=&\Big\{\mathbf{v}_{h}\in (C^{0}(\overline{\Omega }_{f}))^{d}~|~\mathbf{v}_{h|T} \in (\mathbb{P}_1(\mathcal{T}))^d\oplus Span\{\lambda_i^T\}, \forall T\in\mathcal{T}_h,~ \mathbf{v}_{h}|_{\Gamma_f}=0 \Big\},\\
	Q_{fh}&:=&\Big\{q_{h}\in C^0(\overline{\Omega }_{f})~|~q_{h|T}\in\mathbb{P}_1(\mathcal{T}), ~~ \forall T\in\mathcal{T}_h \Big\}, \\
	X_{ph}&:=&\Big\{\psi_{h}\in C^0(\overline{\Omega }_{p})~|~\psi_{h|T}\in\mathbb{P}_1(\mathcal{T}),~\forall T\in\mathcal{T}_h,~ \psi_{h}|_{\Gamma_p}=0 \Big\},
\end{eqnarray*}
where we assume the spaces $\mathbf{X}_{fh}$ and $Q_{fh}$ satisfy the discrete LBB or inf-sup condition \cite{LBB33}, and $\lambda_i^T, i=0,\cdots,d$ are the $d+1$ barycentric coordinate functions for each element $T\in\mathcal{T}_h$. We also define the discrete trace space on the interface $\Gamma$
\begin{eqnarray*}
	Z_{h}&:=&\Big\{\delta_{h}\in C^{0}(\Gamma)~|~\delta_{h|\tau} \in (\mathbb{P}_1(\tau)), \forall \tau\in\mathcal{B}_h,~ \delta_{h}|_{\partial\Gamma}=0 \Big\}.
\end{eqnarray*}

Furthermore, it's easy to find that $Z_{h}$ is the trace space in the sense that
\begin{eqnarray*}
	Y_{fh}:=\mathbf{X}_{fh}|_{\Gamma}\cdot\mathbf{n}_f=Z_{h},\hspace{15mm}Y_{ph}:=X_{ph}|_{\Gamma}=Z_{h}.
\end{eqnarray*}

Similarly to the continuous case, let  $(\phi_{jh},\mathbf{u}_{jh},p_{jh})\in X_{ph}\times \mathbf{X}_{fh}\times Q_{fh}$ are the solution of the finite element approximation of the decoupled Stokes-Darcy system with Robin boundary conditions (\ref{Robinstok1}), (\ref{Robindar1}) and the given parameters $\delta_{pjh}\in Z_{h}$, $\delta_{fjh}\in Z_{h}$. The parallel ensemble domain decomposition finite element method (FE Ensemble DDM) is defined as follows. 

\textbf{\emph{{Parallel FE Ensemble DDM Algorithm}}}

1. Initial values of $\delta_{pjh}^0\in Y_{ph}$, $\delta_{fjh}^0\in Y_{fh}$, $\mathbf{u}^0_{jh}\in \mathbf{X}_{fh}$ and $\phi^0_{jh}\in X_{ph}$ are guessed, which can be zero.

2. For $n=0,1,2,\cdots,$ independently solve the ensemble Stokes and Darcy equations, i.e. $\phi_{jh}^{n+1}\in X_{ph}$ is determined from
\begin{eqnarray}
	\gamma(\overline{\mathbb{K}}\triangledown\phi^{n+1}_{jh},\triangledown\psi_h)_{\Omega_{p}}+g\langle \phi^{n+1}_{jh},\psi_h\rangle
	=\langle\delta^n_{pjh},\psi_h\rangle-\gamma((\mathbb{K}_{j}-\overline{\mathbb{K}})\triangledown\phi^n_{jh},\triangledown\psi_h)_{\Omega_{p}} \hspace{2.5mm} \forall \psi_h \in X_{ph}, \hspace{4mm} \label{FEMdecoupled-3h}
\end{eqnarray}
and $\mathbf{u}_{jh}^{n+1}\in \mathbf{X}_{fh}$ , $p_{jh}^{n+1}\in Q_{fh}$ are determined from
\begin{eqnarray}\label{FEMdecoupled-1h}
	&&a_{f}(\mathbf{u}_{jh}^{n+1},\mathbf{v}_h)+b_f(p_{jh}^{n+1},\mathbf{v}_h)+\gamma\langle\mathbf{u}_{jh}^{n+1}\cdot\mathbf{n}_{f},\mathbf{v}_h\cdot\mathbf{n}_f\rangle
	+\sum_{i=1}^{d-1}\langle\overline{\eta}_{i}\mathbf{u}_{jh}^{n+1}\cdot\tau_{i},\mathbf{v}_h\cdot\tau_{i}\rangle\nonumber\\
	&&=(\mathbf{f}_j,\mathbf{v}_h)_{\Omega_f}-\sum_{i=1}^{d-1}\langle(\eta_{i j}-\overline{\eta}_{i})\mathbf{u}_{jh}^n\cdot\tau_{i},\mathbf{v}_h \cdot\tau_{i}\rangle+\langle\delta_{fjh}^n,\mathbf{v}_h\cdot\mathbf{n}_{f}\rangle \hspace{6.5mm} \forall\mathbf{v}_{h}\in \mathbf{X}_{fh}, \\
	&&b_f(\mathbf{u}_{jh}^{n+1},q_{h})=0 \hspace{73.5mm}\forall q_{h}\in Q_{fh}.\label{FEMdecoupled-2h}
\end{eqnarray}
3. Update $\delta_{fjh}^{n+1} \in Z_{h}$ and $\delta_{pjh}^{n+1} \in Z_{h}$
\begin{eqnarray*}
	&&\delta_{fjh}^{n+1} =a\delta_{pjh}^{n}+bg\phi_{jh}^{n+1},\label{FEMdecoupled-comp1}\\
	&&\delta_{pjh}^{n+1} =c\delta_{fjh}^{n}+d\mathbf{u}_{jh}^{n+1}\cdot\mathbf{n}_f,\label{FEMdecoupled-comp3}
\end{eqnarray*}
where the coefficients
\begin{eqnarray*}
	a=1,\hspace{9mm} b=-2,\hspace{9mm}c=-1,\hspace{9mm} d=2\gamma.
\end{eqnarray*}
Next, we will discuss the convergence of the Parallel FE Ensemble DDM algorithm. Define the following notations for the error functions
\begin{eqnarray*}
	&&\mathbf{e}^n_{ujh}=\mathbf{u}_{jh}-\mathbf{u}^n_{jh},\hspace{9mm}{e}^n_{pjh}=p_{jh}-p^n_{jh},\hspace{10mm}
	{e}^n_{\phi jh}=\phi_{jh}-\phi^n_{jh},\\
	&&{\varepsilon}_{fjh}^{n}=\delta_{fjh}-{\delta}_{fjh}^{n},\hspace{7mm} {\varepsilon}_{pjh}^{n}=\delta_{pjh}-{\delta}_{pjh}^{n},
\end{eqnarray*}
and derive the error equations for any $(\mathbf{v}_{h}, q_{h}, \psi_h)\in (\mathbf{X}_{fh},  Q_{fh},  X_{ph})$
\begin{eqnarray}\label{FEMerr-1} 
	&&a_{f}(\mathbf{e}_{ujh}^{n+1},\mathbf{v}_h)+b_f({e}_{pjh}^{n+1},\mathbf{v}_h)+\gamma\langle\mathbf{e}_{ujh}^{n+1}\cdot\mathbf{n}_{f},\mathbf{v}_h\cdot\mathbf{n}_f\rangle+\sum_{i=1}^{d-1}\langle\overline{\eta}_{i}\mathbf{e}_{ujh}^{n+1}\cdot\tau_{i},\mathbf{v}_h\cdot\tau_{i}\rangle\nonumber\\
	&&=-\sum_{i=1}^{d-1}\langle(\eta_{i j}-\overline{\eta}_{i})\mathbf{e}_{ujh}^n\cdot\tau_{i},\mathbf{v}_h\cdot\tau_{i}\rangle+\langle{\varepsilon}_{fjh}^n,\mathbf{v}_h\cdot\mathbf{n}_{f}\rangle,\\
	&&b_f(\mathbf{e}_{ujh}^{n+1},q_h)=0,\label{FEMerr-2}\\
	&&\gamma(\overline{\mathbb{K}}\triangledown{e}^{n+1}_{\phi jh},\triangledown\psi_h)_{\Omega_{p}}+g\langle {e}^{n+1}_{\phi jh},\psi_h\rangle=\langle{\varepsilon}_{pjh}^n,\psi_h\rangle-\gamma((\mathbb{K}_{j}-\overline{\mathbb{K}})\triangledown{e}^n_{\phi jh},\triangledown\psi_h)_{\Omega_{p}}.\label{FEMerr-3}
\end{eqnarray}
Along the interface $\Gamma$, the error functions can be updated by following formulations
\begin{eqnarray}
	&&\varepsilon_{fjh}^{n+1} =a\varepsilon_{pjh}^{n}+bg{e}_{\phi jh}^{n+1},\label{FEMerrp}\\
	&&\varepsilon_{pjh}^{n+1} =c\varepsilon_{fjh}^{n}+d\mathbf{e}_{ujh}^{n+1}\cdot\mathbf{n}_f.\label{FEMerrf}
\end{eqnarray}
Similar to the continuous case (\ref{ffpp3})-(\ref{ffpp4}), we obtain
\begin{eqnarray}\label{FEMffpp01}
	&&||\varepsilon_{pjh}^{n+1}||^2_{\Gamma}=||\varepsilon_{fjh}^{n}||^2_{\Gamma}-4\gamma a_{f}(\mathbf{e}_{ujh}^{n+1},\mathbf{e}_{ujh}^{n+1})-4\gamma\sum_{i=1}^{d-1}\langle\overline{\eta}_{i}\mathbf{e}_{ujh}^{n+1}\cdot\tau_{i},\mathbf{e}_{ujh}^{n+1}\cdot\tau_{i}\rangle\nonumber\\
	&&\hspace{18mm}-4\gamma\sum_{i=1}^{d-1}\langle(\eta_{i j}-\overline{\eta}_{i})\mathbf{e}_{ujh}^n\cdot\tau_{i},\mathbf{e}_{ujh}^{n+1}\cdot\tau_{i}\rangle,\\
	&&||\varepsilon_{fjh}^{n+1}||^2_{\Gamma}=||\varepsilon_{pjh}^{n}||^2_{\Gamma}-4g\gamma\big[(\overline{\mathbb{K}}\triangledown{e}^{n+1}_{\phi jh},\triangledown{e}_{\phi jh}^{n+1})_{\Omega_{p}}+((\mathbb{K}_{j}-\overline{\mathbb{K}})\triangledown{e}^n_{\phi jh},\triangledown{e}_{\phi jh}^{n+1})_{\Omega_{p}}\big].\label{FEMffpp02}
\end{eqnarray}
\begin{theorem}
	For the case of $\gamma_{f}=\gamma_{p}=\gamma$, with the positive constants $C$, $C_f$, assume that $\overline{\eta}_{i} > \eta_{i}^{\prime \max}$ and $\overline{k}_{\min} > \rho_{\max}^{\prime}$, the Parallel FE Ensemble DDM algorithm has the following convergence
	\begin{eqnarray}\label{FEconvergent}
	&&||\varepsilon_{fjh}^{N}||^2_{\Gamma}+||\varepsilon_{pjh}^{N}||^2_{\Gamma}+||\mathbf{e}_{ujh}^{N}\cdot\tau_{i}||^2_{\Gamma}+||\triangledown{e}^{N}_{\phi jh}||^2_{\Omega_{p}}+||e_{pjh}^N||_{\Omega_{f}}^2+||\triangledown\mathbf{e}_{ujh}^N||^2_{\Omega_{f}}\nonumber\\
&\leq&C\mathcal{E}_h^{N-1}\Big\{||\varepsilon_{fjh}^{0}||^2_{\Gamma}+||\varepsilon_{pjh}^{0}||^2_{\Gamma}+||\mathbf{e}_{ujh}^{0}\cdot\tau_{i}||^2_{\Gamma}+||\triangledown{e}^{0}_{\phi jh}||^2_{\Omega_{p}}\Big\},
		\end{eqnarray}
	where
	\begin{eqnarray}\label{Eh}
		\mathcal{E}_h&=&\max\Big\{1-\frac{2\nu h(2\nu+\frac{C_{f}\gamma}{2\nu}h^{\frac{1}{2}})^{-2}}{1-2\nu h(2\nu+\frac{C_{f}\gamma}{2\nu}h^{\frac{1}{2}})^{-2}}, \hspace{2mm} 1-2\nu h(2\nu+\frac{C_{f}\gamma}{2\nu}h^{\frac{1}{2}})^{-2}, \nonumber\\
		&&\hspace{12mm}\sum_{i=1}^{d-1}\frac{\eta_{i}^{\prime \max}}{(2\overline{\eta}_{i}-\eta_{i}^{\prime \max})}, \hspace{2mm} \frac{\rho_{\max}^{\prime}}{2\overline{k}_{\min}-\rho_{\max}^{\prime}}\Big\}.
	\end{eqnarray}
\end{theorem}
\begin{proof}
	The important point in this proof is an extension operator. 
	 Denote $\mathbf{E}_{fh}\varepsilon_{fjh}^n\in\mathbf{X}_{fh}$ as the discrete Stokes extension of $\varepsilon_{fjh}^n$, and $\mathbf{E}_{fh}\varepsilon_{fjh}^n=\varepsilon_{fjh}^n\cdot\mathbf{n}_f$. According to \cite{Liu34}, we have
	\begin{eqnarray*}
		&&-4\gamma\sum_{i=1}^{d-1}\langle\overline{\eta}_{i}\mathbf{e}_{ujh}^{n+1}\cdot\tau_{i},\mathbf{E}_{fh}\varepsilon_{fjh}^n\cdot\tau_{i}\rangle=0,\hspace*{5mm}b_f({e}_{pjh}^{n+1},\mathbf{E}_{fh}\varepsilon_{fjh}^n)=0,\\
		&&-4\gamma\sum_{i=1}^{d-1}\langle(\eta_{i j}-\overline{\eta}_{i})\mathbf{e}_{ujh}^n\cdot\tau_{i},\mathbf{E}_{fh}\varepsilon_{fjh}^n\cdot\tau_{i}\rangle=0.
	\end{eqnarray*}

Reconsider (\ref{FEMerr-1}), and let $\mathbf{v}_h=\mathbf{E}_{fh}\varepsilon_{fjh}^n$, so we obtain
\begin{eqnarray}\label{FEMerr-11} 
	a_{f}(\mathbf{e}_{ujh}^{n+1},\mathbf{E}_{fh}\varepsilon_{fjh}^n)+\gamma\langle\mathbf{e}_{ujh}^{n+1}\cdot\mathbf{n}_{f},\mathbf{E}_{fh}\varepsilon_{fjh}^n\cdot\mathbf{n}_f\rangle=\langle{\varepsilon}_{fjh}^n,\mathbf{E}_{fh}\varepsilon_{fjh}^n\cdot\mathbf{n}_{f}\rangle.
\end{eqnarray}
Using Cauchy-Schwarz inequality, combining with (\ref{ineq04}) and referring to \cite{Liu34}, we have
\begin{eqnarray}\label{FEMshou} 
	||\varepsilon_{fjh}^n||^2_{\Gamma}
	&=&a_{f}(\mathbf{e}_{ujh}^{n+1},\mathbf{E}_{fh}\varepsilon_{fjh}^n)+\gamma\langle\mathbf{e}_{ujh}^{n+1}\cdot\mathbf{n}_{f},\mathbf{E}_{fh}\varepsilon_{fjh}^n\cdot\mathbf{n}_f\rangle\nonumber\\
	&\leq&a_{f}(\mathbf{e}_{ujh}^{n+1},\mathbf{e}_{ujh}^{n+1})^{\frac{1}{2}}a_{f}(\mathbf{E}_{fh}\varepsilon_{fjh}^n,\mathbf{E}_{fh}\varepsilon_{fjh}^n)^{\frac{1}{2}}+\gamma||\mathbf{e}_{ujh}^{n+1}\cdot\mathbf{n}_{f}||_{\Gamma}||\varepsilon_{fjh}^n||_{\Gamma}\nonumber\\
	&\leq&a_{f}(\mathbf{e}_{ujh}^{n+1},\mathbf{e}_{ujh}^{n+1})^{\frac{1}{2}}(2\nu h^{-\frac{1}{2}})||\varepsilon_{fjh}^n||_{H_{00}^{\frac{1}{2}}(\Gamma)}+\frac{C_{f}\gamma}{2\nu}a_{f}(\mathbf{e}_{ujh}^{n+1},\mathbf{e}_{ujh}^{n+1})^{\frac{1}{2}}||\varepsilon_{fjh}^n||_{\Gamma}\nonumber\\
	&\leq&a_{f}(\mathbf{e}_{ujh}^{n+1},\mathbf{e}_{ujh}^{n+1})^{\frac{1}{2}}(2\nu h^{-\frac{1}{2}})||\varepsilon_{fjh}^n||_{\Gamma}+\frac{C_{f}\gamma}{2\nu}a_{f}(\mathbf{e}_{ujh}^{n+1},\mathbf{e}_{ujh}^{n+1})^{\frac{1}{2}}||\varepsilon_{fjh}^n||_{\Gamma}.
\end{eqnarray}
Then, we have the following fact
\begin{eqnarray}
	||\varepsilon_{fjh}^n||_{\Gamma}\leq(2\nu h^{-\frac{1}{2}})a_{f}(\mathbf{e}_{ujh}^{n+1},\mathbf{e}_{ujh}^{n+1})^{\frac{1}{2}}+\frac{C_{f}\gamma}{2\nu}a_{f}(\mathbf{e}_{ujh}^{n+1},\mathbf{e}_{ujh}^{n+1})^{\frac{1}{2}}.
\end{eqnarray}
Further result can be concluded as follows
\begin{eqnarray}
-4\nu a_{f}(\mathbf{e}_{ujh}^{n+1},\mathbf{e}_{ujh}^{n+1})\leq\big(-4\nu h(2\nu+\frac{C_{f}\gamma}{2\nu}h^{\frac{1}{2}})^{-2}\big)||\varepsilon_{fjh}^n||_{\Gamma}^2.
\end{eqnarray}
Reconsidering (\ref{FEMffpp01}), we deduce 
\begin{eqnarray}\label{FEMerr000}
	&&||\varepsilon_{pjh}^{n+1}||^2_{\Gamma}\leq\big(1-4\nu h(2\nu+\frac{C_{f}\gamma}{2\nu}h^{\frac{1}{2}})^{-2}\big)||\varepsilon_{fjh}^{n}||^2_{\Gamma}-4\gamma\sum_{i=1}^{d-1}\langle\overline{\eta}_{i}\mathbf{e}_{ujh}^{n+1}\cdot\tau_{i},\mathbf{e}_{ujh}^{n+1}\cdot\tau_{i}\rangle\nonumber\\
&&\hspace{18mm}-4\gamma\sum_{i=1}^{d-1}\langle(\eta_{i j}-\overline{\eta}_{i})\mathbf{e}_{ujh}^n\cdot\tau_{i},\mathbf{e}_{ujh}^{n+1}\cdot\tau_{i}\rangle.
\end{eqnarray}
For the term $-4\gamma\sum_{i=1}^{d-1}\langle\overline{\eta}_{i}\mathbf{e}_{ujh}^{n+1}\cdot\tau_{i},\mathbf{e}_{ujh}^{n+1}\cdot\tau_{i}\rangle-4\gamma\sum_{i=1}^{d-1}\langle(\eta_{i j}-\overline{\eta}_{i})\mathbf{e}_{ujh}^n\cdot\tau_{i},\mathbf{e}_{ujh}^{n+1}\cdot\tau_{i}\rangle$, it is similar to the continuous case. We can consider (\ref{cauchy_schwarzf}) to get
\begin{eqnarray}
	&&-4\gamma\sum_{i=1}^{d-1}\langle\overline{\eta}_{i}\mathbf{e}_{ujh}^{n+1}\cdot\tau_{i},\mathbf{e}_{ujh}^{n+1}\cdot\tau_{i}\rangle-4\gamma\sum_{i=1}^{d-1}\langle(\eta_{i j}-\overline{\eta}_{i})\mathbf{e}_{ujh}^n\cdot\tau_{i},\mathbf{e}_{ujh}^{n+1}\cdot\tau_{i}\rangle\nonumber\\
	&&\leq-4\gamma\sum_{i=1}^{d-1}\langle\overline{\eta}_{i}\mathbf{e}_{ujh}^{n+1}\cdot\tau_{i},\mathbf{e}_{ujh}^{n+1}\cdot\tau_{i}\rangle+2\gamma\sum_{i=1}^{d-1}\eta_{i}^{\prime \max}||\mathbf{e}_{ujh}^n\cdot\tau_{i}||^2_{\Gamma}+2\gamma\sum_{i=1}^{d-1}\eta_{i}^{\prime \max}||\mathbf{e}_{ujh}^{n+1}\cdot\tau_{i}||^2_{\Gamma}\nonumber\\
	&&\leq-2\gamma\sum_{i=1}^{d-1}(2\overline{\eta}_{i}-\eta_{i}^{\prime \max})||\mathbf{e}_{ujh}^{n+1}\cdot\tau_{i}||^2_{\Gamma}+2\gamma\sum_{i=1}^{d-1}\eta_{i}^{\prime \max}||\mathbf{e}_{ujh}^n\cdot\tau_{i}||^2_{\Gamma}.
\end{eqnarray}
Let $\mathcal{L}=1-4\nu h(2\nu+\frac{C_{f}\gamma}{2\nu}h^{\frac{1}{2}})^{-2}$. Then, (\ref{FEMerr000}) can be derived
\begin{eqnarray}\label{FEMerr111}
	||\varepsilon_{pjh}^{n+1}||^2_{\Gamma}+2\gamma\sum_{i=1}^{d-1}(2\overline{\eta}_{i}-\eta_{i}^{\prime \max})||\mathbf{e}_{ujh}^{n+1}\cdot\tau_{i}||^2_{\Gamma}
	\leq\mathcal{L}||\varepsilon_{fjh}^{n}||^2_{\Gamma}+2\gamma\sum_{i=1}^{d-1}\eta_{i}^{\prime \max}||\mathbf{e}_{ujh}^n\cdot\tau_{i}||^2_{\Gamma}.\hspace{2mm}
\end{eqnarray}
Meanwhile, considering (\ref{FEMffpp02}) and (\ref{cauchy_schwarzp}), we have
\begin{eqnarray}\label{FEMerr222}
	||\varepsilon_{fjh}^{n+1}||^2_{\Gamma}+2g\gamma(2\overline{k}_{\min}-\rho_{\max}^{\prime})||\triangledown{e}^{n+1}_{\phi jh}||^2_{\Omega_{p}}\leq||\varepsilon_{pjh}^{n}||^2_{\Gamma}+2g\gamma\rho_{\max}^{\prime}||\triangledown{e}^{n}_{\phi jh}||^2_{\Omega_{p}}.
\end{eqnarray}
Multiply $\frac{1+\mathcal{L}}{2}$ with (\ref{FEMerr222}) and combine with (\ref{FEMerr111}) to derive the following fact
\begin{eqnarray}
	&&\frac{1+\mathcal{L}}{2}||\varepsilon_{fjh}^{n+1}||^2_{\Gamma}+||\varepsilon_{pjh}^{n+1}||^2_{\Gamma}+2\gamma\sum_{i=1}^{d-1}(2\overline{\eta}_{i}-\eta_{i}^{\prime \max})||\mathbf{e}_{ujh}^{n+1}\cdot\tau_{i}||^2_{\Gamma}\nonumber\\
	&&+\frac{1+\mathcal{L}}{2}2g\gamma(2\overline{k}_{\min}-\rho_{\max}^{\prime})||\triangledown{e}^{n+1}_{\phi jh}||^2_{\Omega_{p}}\nonumber\\
	&\leq&\mathcal{L}||\varepsilon_{fjh}^{n}||^2_{\Gamma}+\frac{1+\mathcal{L}}{2}||\varepsilon_{pjh}^{n}||^2_{\Gamma}+2\gamma\sum_{i=1}^{d-1}\eta_{i}^{\prime \max}||\mathbf{e}_{ujh}^n\cdot\tau_{i}||^2_{\Gamma}\nonumber\\
	&&+\frac{1+\mathcal{L}}{2}2g\gamma\rho_{\max}^{\prime}||\triangledown{e}^{n}_{\phi jh}||^2_{\Omega_{p}},
\end{eqnarray}
According to Lemma \ref{a1b1c1<a2b2c2}, it is easy to check
\begin{eqnarray*}
	&&a_{2h}=\mathcal{L}<a_{1h}=\frac{1+\mathcal{L}}{2},\hspace{5mm}c_{2h}=2\gamma\sum_{i=1}^{d-1}\eta_{i}^{\prime \max}<c_{1h}=2\gamma\sum_{i=1}^{d-1}(2\overline{\eta}_{i}-\eta_{i}^{\prime \max}), \\
	&&b_{2h}=\frac{1+\mathcal{L}}{2}<b_{1h}=1,\hspace{6mm} d_{2h}=\frac{1+\mathcal{L}}{2}2g\gamma\rho_{\max}^{\prime}<d_{1h}=\frac{1+\mathcal{L}}{2}2g\gamma(2\overline{k}_{\min}-\rho_{\max}^{\prime}).
\end{eqnarray*} 
Hence, we arrive at
\begin{eqnarray}
	&&\frac{1+\mathcal{L}}{2}||\varepsilon_{fjh}^{N}||^2_{\Gamma}+||\varepsilon_{pjh}^{N}||^2_{\Gamma}\nonumber\\
	&&+2\gamma\sum_{i=1}^{d-1}(2\overline{\eta}_{i}-\eta_{i}^{\prime \max})||\mathbf{e}_{ujh}^{N}\cdot\tau_{i}||^2_{\Gamma}+\frac{1+\mathcal{L}}{2}2g\gamma(2\overline{k}_{\min}-\rho_{\max}^{\prime})||\triangledown{e}^{N}_{\phi jh}||^2_{\Omega_{p}}\nonumber\\
	&\leq&\mathcal{E}_h^{N-1}\Big\{\mathcal{L}||\varepsilon_{fjh}^{0}||^2_{\Gamma}+\frac{1+\mathcal{L}}{2}||\varepsilon_{pjh}^{0}||^2_{\Gamma}\nonumber\\
	&&+2\gamma\sum_{i=1}^{d-1}\eta_{i}^{\prime \max}||\mathbf{e}_{ujh}^{0}\cdot\tau_{i}||^2_{\Gamma}+\frac{1+\mathcal{L}}{2}2g\gamma\rho_{\max}^{\prime}||\triangledown{e}^{0}_{\phi jh}||^2_{\Omega_{p}}\Big\},
\end{eqnarray}
where $\mathcal{E}_h$ has been Adefined in (\ref{Eh}).

We can further analysis the error equation (\ref{FEMerr-1}) and  (\ref{FEMerr-2})  to derive that the  velocity $||\triangledown\mathbf{e}_{ujh}^n||_{\Omega_{f}}^2$ can be controlled by $||\varepsilon_{fjh}^{n}||^2_{\Gamma}$ and $||\mathbf{e}_{ujh}^{n}\cdot\tau_{i}||^2_{\Gamma}$.
Subsequently, similar to (\ref{pshoulian}), we can get the convergence of $||e_{pjh}^n||^2_{\Omega_{f}}$. Overall, we can summarize the results of the proposed algorithm convergence (\ref{FEconvergent}).
\end{proof}



 As for the case $\gamma_{f}<\gamma_{p}$, since the convergence rate is mesh-independent, the result of finite element form can directly refer to Theorem \ref{case2<thm}.

\begin{rem}
	This remark will discuss the value of $\mathcal{E}_h$ to demonstrate that the convergence of our algorithm is related to the mesh size in the case $\gamma_{f}=\gamma_{p}=\gamma$. From above analysis, we need to assume  $\overline{\eta}_{i} > \eta_{i}^{\prime \max}$ and $\overline{k}_{\min} > \rho_{\max}^{\prime}$, that means a small perturbation for $\mathbb{K}$. 
	In other words, the value of $\eta_{i}^{\prime \max}$ and $\rho_{\max}^{\prime}$ need to be small enough. Thus, the last two terms $\frac{\sum_{i=1}^{d-1}\eta_{i}^{\prime \max}}{\sum_{i=1}^{d-1}(2\overline{\eta}_{i}-\eta_{i}^{\prime \max})}$ and $\frac{\rho_{\max}^{\prime}}{2\overline{k}_{\min}-\rho_{\max}^{\prime}}$ must not be the maximum of all the terms. For the first two terms $1-\frac{2\nu h(2\nu+\frac{C_{f}\gamma}{2\nu}h^{\frac{1}{2}})^{-2}}{1-2\nu h(2\nu+\frac{C_{f}\gamma}{2\nu}h^{\frac{1}{2}})^{-2}}$ and $1-2\nu h(2\nu+\frac{C_{f}\gamma}{2\nu}h^{\frac{1}{2}})^{-2}$, as $h$ decreases, through simple algebraic operations, it can be discovered that they both approach to $1$. In summary, regardless of which term $\mathcal{E}_h$ is taken, the convergence rate of the parallel FE ensemble DDM algorithm is dependent on the mesh size.
	
	 Based on above discussion, we can provide an guidance to select Robin parameters $\gamma_f$ and $\gamma_p$. Since the convergence rate of case $\gamma_p=\gamma_f$ depends on the mesh, while the convergence rate of case $\gamma_f < \gamma_p$ is mesh-independent, we suggest the following: to enhance efficiency, use the Robin parameters selection scheme (\ref{opm}) in Theorem 4.4 for numerical experiments.
\end{rem}

	\section{ Multi-level Monte Carlo Ensemble DDM}
	As stated above, we focus on proposing and analyzing the ensemble DDM algorithm after applying the Monte Carlo method to the random PDEs. The classical Monte Carlo (MC) method converges slowly and requires many samples $J$. This makes the fast solution of the Stokes-Darcy coupled problems with random hydraulic conductivity and body force unrealistic. Therefore, we adopt the multi-level Monte Carlo (MLMC) method for the probability space in this section, which results in a more efficient numerical solution with lower computational effort than MC.  We also  present data evidence in the numerical experiments part.
	

	To fit in the hierarchic nature of MLMC ensemble DDM, we can consider the same regular and quasi-uniform trangulations $\mathcal{T}_{fh}^l$ and $\mathcal{T}_{ph}^l$ for the subdomains $\Omega_f, \ \Omega_p$ at different mesh level $l$ ($l=0,1,\cdots,L$, where $L$ is the largest level of the selected finest grid). Similarly, we can define the corresponding FE space at different mesh level $l$ as $\mathbf{X}_{fh}^l, Q_{fh}^l$ and $X_{ph}^l$. Note that the sequence of FE spaces satisfies
	\begin{eqnarray*}
		&&\mathbf{X}_{fh}^0 \subset \mathbf{X}_{fh}^1 \subset \cdots \subset \mathbf{X}_{fh}^l \subset \cdots \subset \mathbf{X}_{fh}^L,\\
		&&{Q}_{fh}^0 \subset {Q}_{fh}^1 \subset \cdots \subset {Q}_{fh}^l \subset \cdots \subset {Q}_{fh}^L,\\
		&&{X}_{ph}^0 \subset {X}_{ph}^1 \subset \cdots \subset {X}_{ph}^l \subset \cdots \subset {X}_{ph}^L.
	\end{eqnarray*}
	Let $\mathbf{u}_{fh}^l (\mathbf{x},\omega) \in \mathbf{X}_{fh}^l$ be the numerical velocity of Stokes at each mesh level $l$. We can write the numerical pressure of Stokes and the piezometric head of Darcy as  ${p}_{fh}^l (\mathbf{x},\omega) \in {Q}_{fh}^l$ and ${\phi}_{ph}^l (\mathbf{x},\omega) \in {X}_{ph}^l$, respectively. To simplify, we only use the Stokes velocity to explain the basic idea of MLMC, and the same methods apply to other quantities. The Stokes velocity of MLMC ensemble DDM on the $L$-th level grid is given by
	\begin{eqnarray*}
		\mathbf{u}_{fh}^L(\mathbf{x},\omega)=\sum_{l=1}^{L}(\mathbf{u}_{fh}^{l}(\mathbf{x},\omega)-\mathbf{u}_{fh}^{l-1}(\mathbf{x},\omega))+\mathbf{u}_{fh}^{0}(\mathbf{x},\omega).
	\end{eqnarray*}
	By linearity of the expectation operator $E_{J^l}[\cdot]$, where $J^l$ is the number of samples, we have
	\begin{eqnarray}\label{MLMC_SV}
		E_{J^L}[\mathbf{u}_{fh}^L(\mathbf{x},\omega)]&=&E_{J^L}\Big[\sum_{l=1}^{L}(\mathbf{u}_{fh}^{l}(\mathbf{x},\omega)-\mathbf{u}_{fh}^{l-1}(\mathbf{x},\omega))+\mathbf{u}_{fh}^{0}(\mathbf{x},\omega)\Big]\nonumber\\
		&=&\sum_{l=1}^{L}E_{J^l}\Big[(\mathbf{u}_{fh}^{l}(\mathbf{x},\omega)-\mathbf{u}_{fh}^{l-1}(\mathbf{x},\omega))\Big]+E_{J^0}[\mathbf{u}_{fh}^{0}(\mathbf{x},\omega)].
	\end{eqnarray}
	where $	E_{J^l}[\mathbf{u}_{fh}^l(\mathbf{x},\omega)]\approx\frac{1}{J^l}\sum^{J^l}_{j=1}\mathbf{u}_{fh}^l(\mathbf{x},\omega_j)$. MLMC typically has a smaller number of samples in the $l$-th level grid than in the $(l-1)$-th level grid, i.e. $J^l < J^{l-1}$. Therefore, to estimate $E_{J^l}\Big[(\mathbf{u}_{fh}^{l}(\mathbf{x},\omega)-\mathbf{u}_{fh}^{l-1}(\mathbf{x},\omega))\Big]$, we usually randomly select $J^l$ solutions from the  $(l-1)$-th level grid samples. 
	To obtain the expected numerical solutions on the $L$-level grid, we assign the number of samples $J^l$ for every $l$-level grid, and apply our ensemble DDM algorithm from Section 3 to each grid level. As the grid level rises, the number of samples we assign drops exponentially, resulting in significant computational savings.
	
	\textbf{\emph{{MC FE Ensemble DDM Algorithm}}}
	
		1. Generate a number of independently, identically distributed (i.i.d) samples for the random hydraulic conductivity $\mathbb{K}(\mathbf{x},\omega_j)$ and the random forces $\mathbf{f}(\mathbf{x},\omega_j)$, where $j=1,\cdots,J$.
	
		2. Apply the Parallel FE Ensemble DDM algorithm to solve for approximate solutions $\mathbf{u}_{fh}(\mathbf{x},\omega_j)$, $p_{fh}(\mathbf{x},\omega_j)$, $\psi_{ph}(\mathbf{x},\omega_j), j=1,\cdots,J$.
	
		3. Output required statistical information, such as the expectation of the Stokes velocity $\mathbf{u}_{fh}(\mathbf{x},\omega): E_{J}[\mathbf{u}_{fh}(\mathbf{x},\omega)]\approx\frac{1}{J}\sum^{J}_{j=1}\mathbf{u}_{fh}(\mathbf{x},\omega_j)$.

\textbf{\emph{{MLMC FE Ensemble DDM Algorithm}}}
	
		1. Generate a number of independently, identically distributed (i.i.d) samples for the random hydraulic conductivity $\mathbb{K}^l(\mathbf{x},\omega_j)$ and the random forces $\mathbf{f}^l(\mathbf{x},\omega_j)$ at every different mesh level $l$, where $j=1,\cdots,J^l, l=1,\cdots,L$, and $L$ is the finest grid level. Typically, $J^l\leq J^{l-1}$.
	
		2. Apply the Parallel FE Ensemble DDM algorithm on every $l$ level grid to solve for approximate solutions $\mathbf{u}_{fh}^l(\mathbf{x},\omega_j), p_{fh}^l(\mathbf{x},\omega_j), \psi_{ph}^l(\mathbf{x},\omega_j), j=1,\cdots,J^l, l=0,\cdots,L$.
	
		3. Output required statistical information for the Stokes velocity is given by (\ref{MLMC_SV}). Applying the same methods to other quantities.
	
	To show the benefits of MLMC in numerical analysis and the influence of key parameters like sample number $J^l$ and mesh size $h^l$ on convergence rate, we use the following theorem, which mainly follows the analysis technique of \cite{Mul02} (Lemma 4.4).
	\begin{theorem}
		Let $(\mathbf{u}_{f}, p_f; \phi_p)$ denote the solution of the random Stokes-Darcy problems. Suppose that $(\check{\mathbf{u}}_{fh}^l, \check{p}_{fh}^l; \check{\phi}_{ph}^l)$ and $(\tilde{\mathbf{u}}_{fh}^l, \tilde{p}_{fh}^l; \tilde{\phi}_{ph}^l)$ are the approximations of expect value by MC ensemble DDM and MLMC ensemble DDM, respectively. Then, the error bounds for random Stokes-Darcy problems are given as follows:
		\begin{eqnarray}
			||\nabla \mathbb{E}\left[\mathbf{u}_{f}\right]- \nabla \check{\mathbf{u}}_{fh}^{l}||_{\Omega_{f}}+||\mathbb{E}\left[p_f\right]-\check{p}_{fh}^{l}||_{\Omega_f} + ||\nabla \mathbb{E}\left[{\phi}_{p}\right]- \nabla \check{\phi}_{ph}^{l}||_{\Omega_{p}}&\leq& C \big(h^l+(J^l)^{-\frac{1}{2}}\big),\nonumber \\	
			||\nabla \mathbb{E}\left[\mathbf{u}_{f}\right]-\nabla \tilde{\mathbf{u}}_{fh}^{l}||_{\Omega_{f}}+||\mathbb{E}\left[p_f\right]-\tilde{p}_{fh}^{l}||_{\Omega_f} + ||\nabla \mathbb{E}\left[{\phi}_{p}\right]-\nabla \tilde{\phi}_{ph}^{l}||_{\Omega_{p}}&\leq& C\big(h^L+\sum_{l=0}^{L}h^l(J^l)^{-\frac{1}{2}}\big),\nonumber
		\end{eqnarray}
		where $C$ depends on $\mathbf{u}_{f}, p_f$ and $ \phi_p$.
	\end{theorem}
	\begin{rem}
		We can derive the above theorem in a similar manner to the proof of Lemma 4.4 in \cite{Mul02}, but by utilizing the linearity of $\nabla \mathbb{E}\left[\mathbf{u}_{f}\right]=\mathbb{E}\left[\nabla\mathbf{u}_{f}\right]$. Because the true solutions of Stokes velocity and Darcy piezometric head have different regularity in our paper and in Proposition 4.3 of \cite{Mul02}, we have slightly adjusted the above theorem. Moreover, \cite{Mul01} presents the MLMC ensemble algorithm for random parabolic PDEs. Without the time correlation in Theorem 5 of \cite{Mul01}, we can reach a similar result as above.
	\end{rem}

\section{Numerical Experiments}


In this section, we present the results of three numerical experiments to demonstrate the accuracy and efficiency of the proposed ensemble DDM. The first experiment involves simulating a problem with a known exact solution to check the convergence of the algorithm and demonstrate its optimal convergence accuracy. The second example presented aims to demonstrate the effectiveness and efficiency of the ensemble algorithm combined with the Monte Carlo method in simulating the Stokes-Darcy system with a random hydraulic conductivity tensor. We also confirm that MLMC ensemble DDM offers more efficient numerical solutions based on the overall computational effort than MC. Moreover, we simulate the coupling of the “3D Shallow Water” system with random hydraulic conductivity to present the complicated flow characteristics. This is a common challenge in modeling porous media flow, where the hydraulic conductivity tensor is often an uncertain parameter that can greatly impact the accuracy of the model’s predictions. 

All numerical tests were carried out by using the open-source finite element method software FreeFEM++ \cite{F18}.

\subsection{Smooth Problem with Convergence Test}
The first example with exact solutions is adapted from \cite{JiangN29} to verify the convergence and to check the feasibility of Ensemble DDM. The free fluid flow domain $\Omega_{f}= [0,\pi]\times[0,1]$ and the porous medium domain $\Omega_{p}=[0,\pi]\times[-1,0]$ are considered, including the interface  $\Gamma= \{ 0 \leq x \leq \pi, y=0\}$. For the computational convenience, we assume $z=0$, and  other physical parameters $ \nu, \gamma$ and $\alpha$ to be $ 1.0 $. Moreover, we can follow Theorem 4.4 to
choose the optimized Robin parameters $\gamma_{f}^{\ast},\ \gamma_{p}^{\ast}$. The random hydraulic conductivity tensor $\mathbb{K}$ will be assumed as
\begin{eqnarray*}
	\mathbb{K}=\mathbb{K}_{j}=\left[\begin{array}{cc}
		k_{11}^{j} & 0 \\
		0 & k_{22}^{j}
	\end{array}\right], \quad j=1, \ldots, J,
\end{eqnarray*}
where $\mathbb{K}_{j}$ is one of the samples of the random hydraulic conductivity tensor $\mathbb{K}$.
The exact solution is selected as
\begin{eqnarray*}
	\begin{aligned}
		\phi_{p} &=\left(e^{y}-e^{-y}\right) \sin (x),\hspace{6mm}
			p_{f} =0,\\
		\mathbf{u}_{f} &=\left[\frac{k_{11}^{j}}{\pi} \sin (2 \pi y) \cos (x),\Big(-2 k_{22}^{j}+\frac{k_{22}^{j}}{\pi^{2}} \sin ^{2}(\pi y)\Big) \sin (x)\right]^{T}.
	\end{aligned}
\end{eqnarray*}

For better comparison with \cite{JiangN29}, we can choose the same group of simulations with $J=3$ members, which correspond to different hydraulic conductivity, i.e. $k_{11}^{1}=k_{22}^{1}=2.21,~ k_{11}^{2}=k_{22}^{2}=4.11,~ k_{11}^{3}=k_{22}^{3}=6.21$. In Table \ref{OrderDDME1}, we present numerical errors of the velocity in $L^2$-norm, $H^1$-norm, and that of the pressure in $L^2$-norm for free fluid flow. As for the porous medium flow, the piezometric head in $L^2$-norm, $H^1$-norm are also shown. From Table \ref{OrderDDME1}, we can see that our algorithm has $h$-independent convergence rate in the case of $\gamma_f<\gamma_p$, which supports Theorem \ref{case2<thm}. Numerical results in Table \ref{OrderDDME1} also demonstrate the optimal convergence for both velocity and pressure, which has the same error accuracy as 
\cite{JiangN29}.

	\begin{table}[!h]
		\caption{Convergence performance for Ensemble DDM ($P1b-P1-P1$ and $\gamma_{f}=\gamma_{f}^{\ast},\gamma_{p}=\gamma_{p}^{\ast}$).}
		\label{OrderDDME1}\tabcolsep 0pt \vspace*{-20pt}
		\par
		\begin{center}
			\def\temptablewidth{1.0\textwidth}
			{\rule{\temptablewidth}{1pt}}
			\begin{tabular*}{\temptablewidth}{@{\extracolsep{\fill}}c|ccccccc}
				\hline
				$k_j$&$h $& itr & $\frac{||\mathbf u_f-\mathbf u_{jh}||_{0}}{||\mathbf u_f||_{0}}$ & $\frac{||\mathbf u_f-\mathbf u_{jh}||_{1}}{||\mathbf u_f||_{1}}$ & $\frac{||p-p_{jh} ||_{0}}{||p||_0}$  &  $\frac{||\phi-\phi_{jh}||_{0}}{||\phi||_0}$   & $\frac{||\phi-\phi_{jh}||_{1}}{||\phi||_1}$ \\
				\hline
				2.21 & $\frac{1}{16} $&   21  &  0.011867 &  0.178799  &  0.358833  & 0.007028  &  0.079797
				\\ 
				& $\frac{1}{32}$  & 21    &   0.003165  &  0.091732 & 0.140762  &      0.001980  &    0.041945 
				\\ 
				& $\frac{1}{64}$ & 21      &  0.000752  &   0.044213   &      0.057421  &  0.000451 &   0.020364
				\\
				& $\frac{1}{128}$ &21     & 0.000186   &  0.021963 &  0.027714  & 0.000114 &  0.010205 
				\\
				\hline
				4.11 & $\frac{1}{16} $& 7    & 0.011852  & 0.178807 & 0.358478     &   0.007047  & 0.079798
				\\ 
				& $\frac{1}{32}$  & 7    &  0.003162 & 0.091733 & 0.140689 &    0.001985 &  0.041945
				\\ 
				& $\frac{1}{64}$ &7     & 0.000752  &  0.044213 & 0.057413 &    0.000451 &  0.020364 
				\\
				& $\frac{1}{128}$ & 7     & 0.000185 &   0.021962 & 0.027712  &      0.000113 &   0.010205
				\\
				\hline
				6.21 & $\frac{1}{16} $& 20    & 0.011843  &  0.178811  & 0.358368  & 0.007052 & 0.079798 
				\\ 
				& $\frac{1}{32}$  &20    & 0.003160 & 0.091734&  0.140666  &    0.001986 & 0.041945
				\\ 
				& $\frac{1}{64}$ &20    &  0.000751 & 0.044213 &  0.057410 &    0.000451 & 0.020364 
				\\
				& $\frac{1}{128}$ &   20   & 0.000185 & 0.021962  & 0.027712   &    0.000113  &  0.010205
				\\
				\hline
			\end{tabular*}%
		\end{center}
	\end{table}

\begin{figure}[htbp]\label{test1}
	\centering
	\subfigure{
		\centering
		\includegraphics[width=0.3\textwidth]{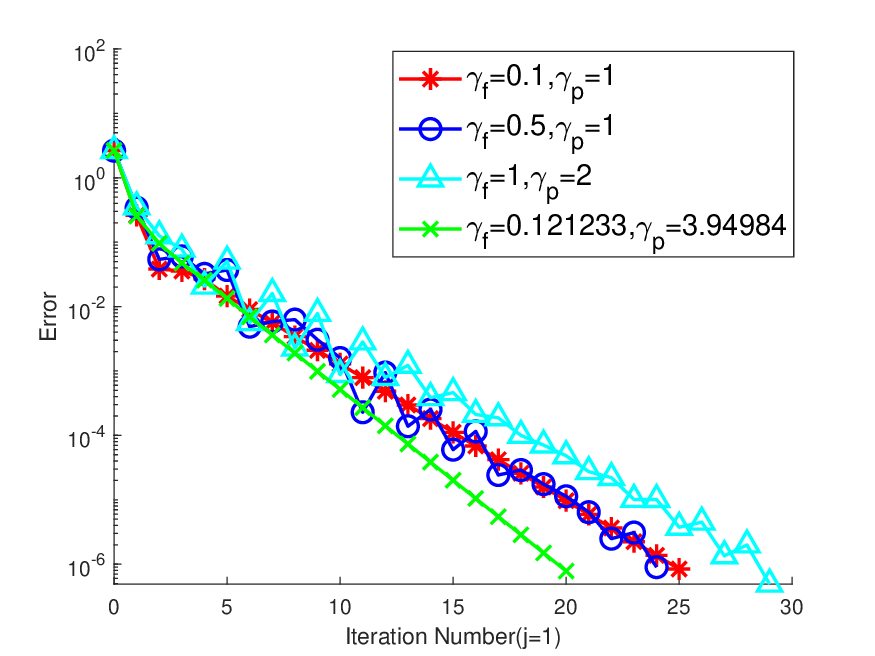}
}
\subfigure{
	\centering
	\includegraphics[width=0.3\textwidth]{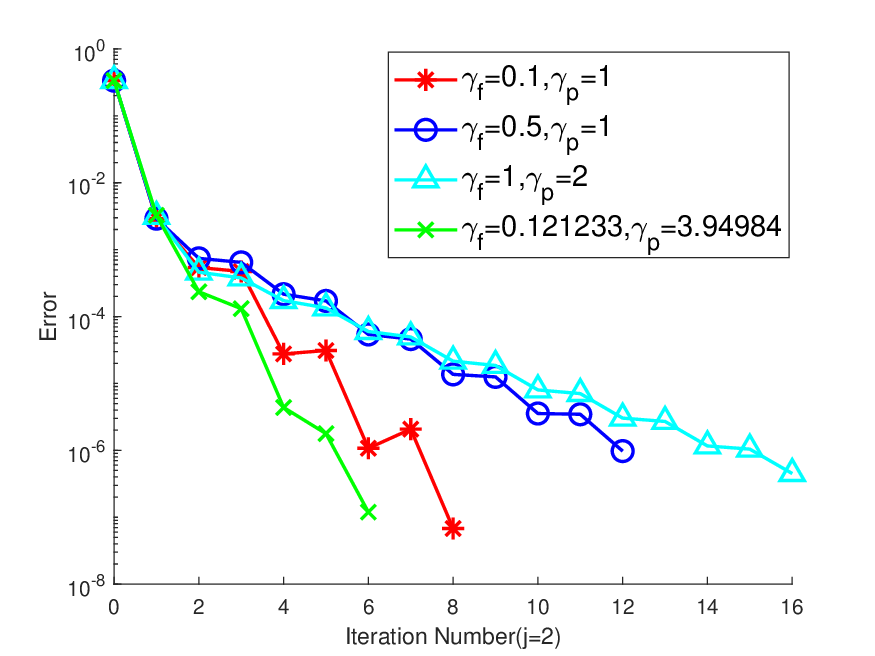}
}
\subfigure{
	\centering
	\includegraphics[width=0.3\textwidth]{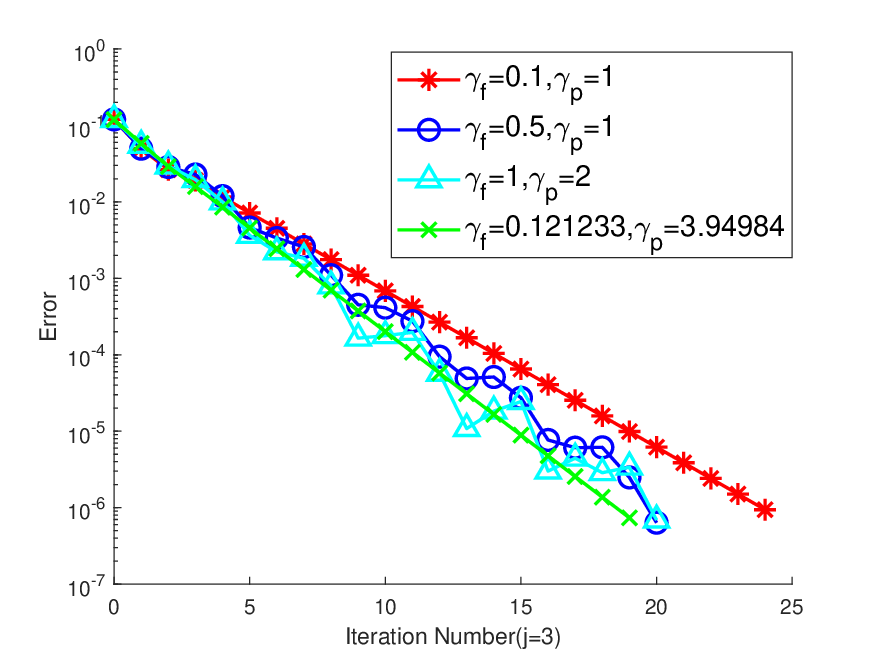}
}
\centering
\caption{The iterates of the Ensemble DDM (J=3) with different Robin parameters $\gamma_{f}$ and $\gamma_{p}$ while $h=\frac{1}{32}$.}
\end{figure}

In addition, to verify the theoretical analysis regarding the performance of the optimal Robin parameters, we further test the proposed Ensemble DDM with different Robin parameters $\gamma_f$ and $\gamma_p$ while $h=\frac{1}{32}$, and  show the comparison of the corresponding numbers of the iterations. We exhibit the evolution of the $L^2$-error about the successive solution components of the Stokes velocities with increasing iteration steps in Fig. \ref{test1}. We choose the optimal Robin parameters according to Theorem \ref{remark} and the other four groups of the Robin parameters to carry out numerical experiments. Among the compared Robin parameters, the optimal one will converge fast and have fewer iterations, which supports our theoretical analysis.

The proof of Theorem 4.4 has strong constraints and the above example happens to meet them all. Hence, we next want to check if the constraints are necessary in numerical experiments. We choose the group of simulations with $J=3$ members, which correspond to different hydraulic conductivity, i.e. $k_{11}^1=2.11, k_{22}^1=3.11$, $k_{11}^2=4.11, k_{22}^2=5.21$, $k_{11}^3=6.21, k_{22}^3=1.21$. It is clearly to see that $k_{11}^j\neq k_{22}^j$ and the sink/source term ${f}_{pj}\neq0$, where $j=1,2,3$. From Fig. 7.2, the optimal Robin parameters case can also converge fast and have fewer iterations.
\begin{figure}[htbp]\label{k1neqk2}
	\centering
	\subfigure{
			\centering
			\includegraphics[width=0.3\textwidth]{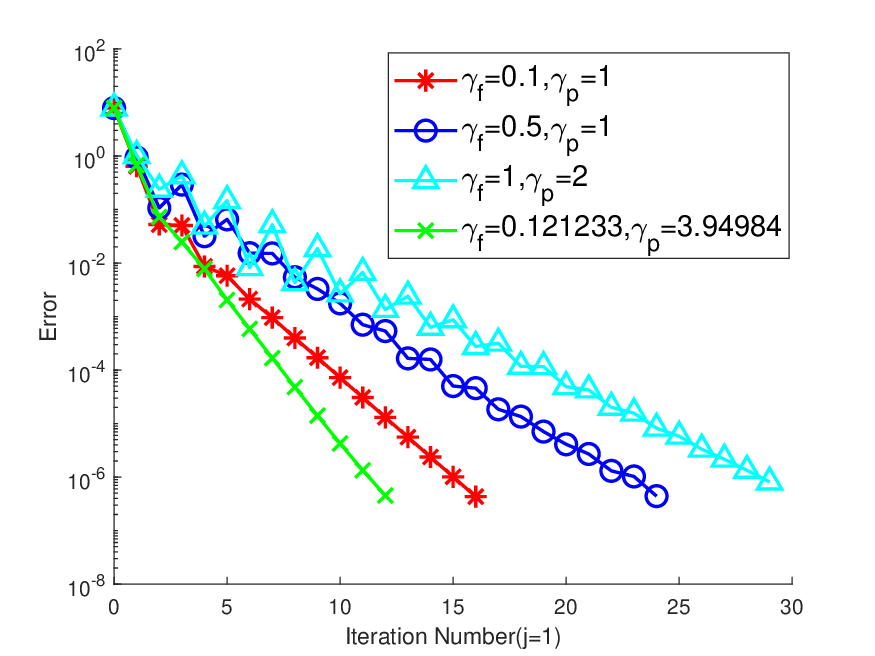}
}
	\subfigure{
			\centering
			\includegraphics[width=0.3\textwidth]{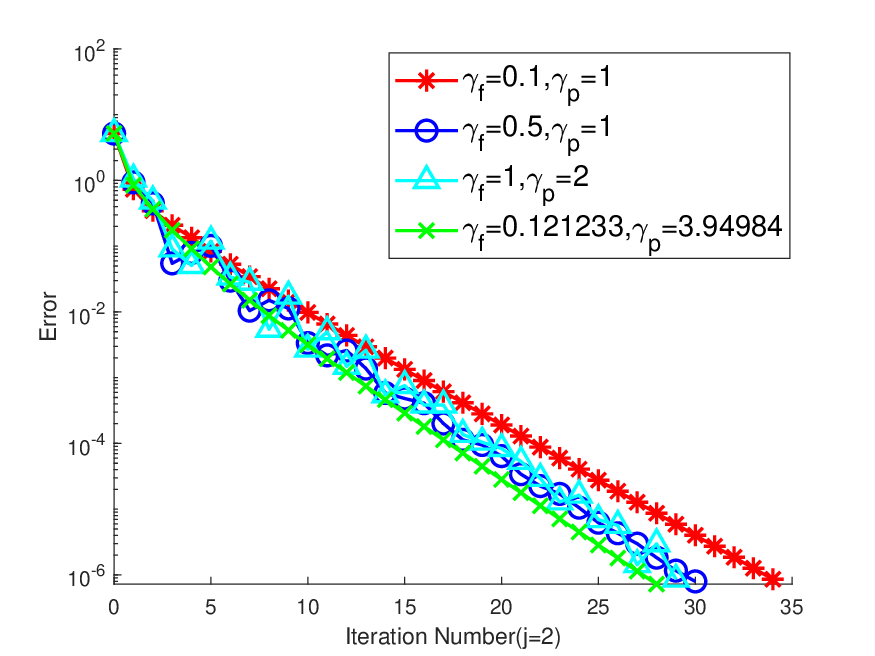}
}
	\subfigure{
			\centering
			\includegraphics[width=0.3\textwidth]{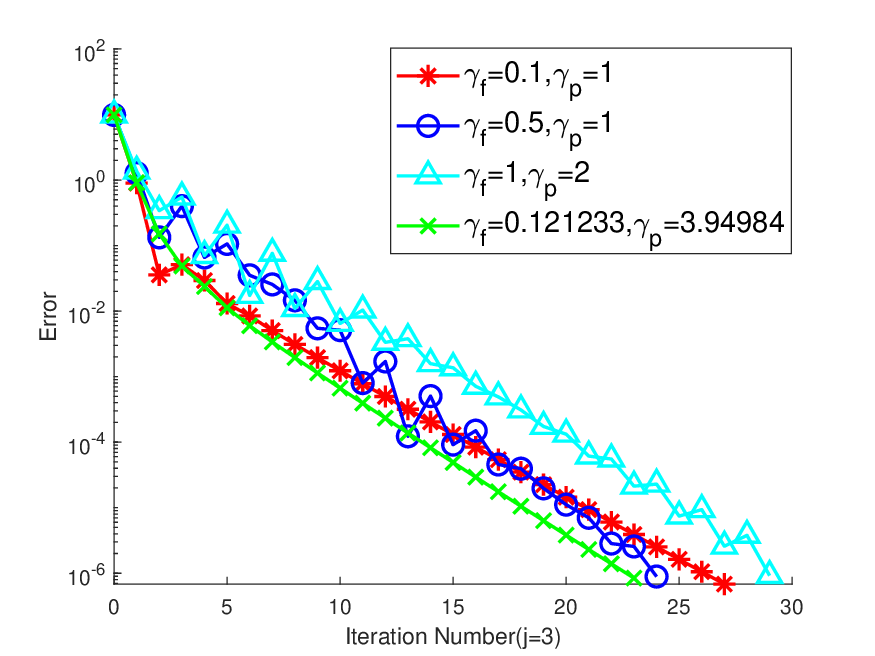}
}
	\centering
	\caption{The iterates of the Ensemble DDM (J=3) with different Robin parameters $\gamma_{f}$ and $\gamma_{p}$ while $h=\frac{1}{32}$.}
\end{figure}
\subsection{Model Problem with Random Hydraulic Conductivity Tensor}
We further consider one model problem with a random hydraulic conductivity tensor that depends on spatial coordinates. To deal with the random model, we can use the Monte Carlo method for sampling and the proposed ensemble DDM for numerical simulations. The random hydraulic conductivity tensor $\mathbb{K}$ can be constructed as follows, which varies in the vertical direction.
\begin{eqnarray*}
	\begin{aligned}
		&\mathbb{K}(\mathbf{x}, \omega)=\left[\begin{array}{cc}
			k_{11}(\mathbf{x}, \omega) & 0 \\
			0 & k_{22}(\mathbf{x}, \omega)
		\end{array}\right], \quad \text { and } \\
		&k_{11}(\mathbf{x}, \omega)=k_{22}(\mathbf{x}, \omega)=a_{0}+\sigma \sqrt{\lambda_{0}} Y_{0}(\omega)+\sum_{i=1}^{n_{f}} \sigma \sqrt{\lambda_{i}}\left[Y_{i}(\omega) \cos (i \pi y)+Y_{n_{f}+i}(\omega) \sin (i \pi y)\right],
	\end{aligned}
\end{eqnarray*}
where $\mathbf{x}=(x, y)^{T}, \lambda_{0}=\frac{\sqrt{\pi L_{c}}}{2}, \lambda_{i}=\sqrt{\pi} L_{c} e^{-\frac{\left(i \pi L_{c}\right)^{2}}{4}}$ for $i=1, \ldots, n_{f}$ and $Y_{0}, \ldots, Y_{2 n_{f}}$ are uncorrelated random variables with unit variance and zero means. In the following numerical test, we can choose the desired physical correlation length $L_{c}=0.25$ for the random field and $a_{0}=1, \sigma=0.15, n_{f}=3$. The random variables $Y_{0}, \ldots, Y_{2 n_{f}}$ are assumed independent and uniformly distributed in the interval $[-\sqrt{3}, \sqrt{3}] $. With the above assumptions, we note that the random functions $k_{11}(\vec{x}, \omega)$ and $k_{22}(\vec{x}, \omega)$ are positive, and the corresponding random hydraulic conductivity tensor $\mathbb{K}$ is symmetric positive definite.
The domain and parameters are the same as those in the first test. However, in this test, the forcing terms are defined as follows:
\begin{eqnarray*}
	\begin{aligned}
	f_p &=\left(e^{y}-e^{-y}\right) \sin (x),\\
	f_{f1}&=[(1+\nu+4\nu\pi^2)\frac{k_{11}(\mathbf{x}, \omega)}{\pi}]sin (2\pi y)cos(x),\\
	f_{f2}&=-2\nu k_{22}(\mathbf{x}, \omega)cos (2\pi y)sin (x)+(1+\nu)[-2k_{22}(\mathbf{x}, \omega)+\frac{k_{22}(\mathbf{x}, \omega)}{\pi^2}sin^2 (\pi y)]sin (x).
\end{aligned}
\end{eqnarray*}
The Dirichlet boundary conditions need to satisfy
\begin{eqnarray*}
	\begin{aligned}
		\phi_p &=\left(e^{y}-e^{-y}\right) \sin (x),\\
		{\mathbf{u}_f}&=\Big[\frac{k_{11}(\mathbf{x}, \omega)}{\pi}sin (2\pi y)cos (x), \Big(-2k_{22}(\mathbf{x}, \omega)+\frac{k_{22}(\mathbf{x}, \omega)}{\pi^2}sin^2(\pi y)\Big)sin (x)\Big]^T.
	\end{aligned}
\end{eqnarray*}

We first check the rate of convergence with respect to the number of samples $J$. As the exact solution to the random Stokes-Darcy system is unknown, we take the ensemble mean of numerical solutions of $J_0=1000$ realizations as our exact solution (expectation) $(\mathbf{u}_{f0}, p_{f0}, \phi_{p0})$ and then evaluate the approximation errors based on this. The numerical results with $J=40,60,80,140,220$ realizations are tested, and the errors of $||\mathbf{u}_{f0}-\mathbf{u}_{fh}||_0, ||\mathbf{u}_{f0}-\mathbf{u}_{fh}||_1, ||p_{f0}-p_{fh}||_0, ||\phi_{p0}-\phi_{ph}||_0, ||\phi_{p0}-\phi_{ph}||_1$ are plotted in Fig. 7.3 by using linear regression. It is seen that the rate of convergence with respect to $J$ is close to 0.5, which coincides with the  Monte Carlo convergence results in Theorem 6.1.
\begin{figure}[htbp]\label{errtest}
	\centering
		\includegraphics[width=0.3\textwidth]{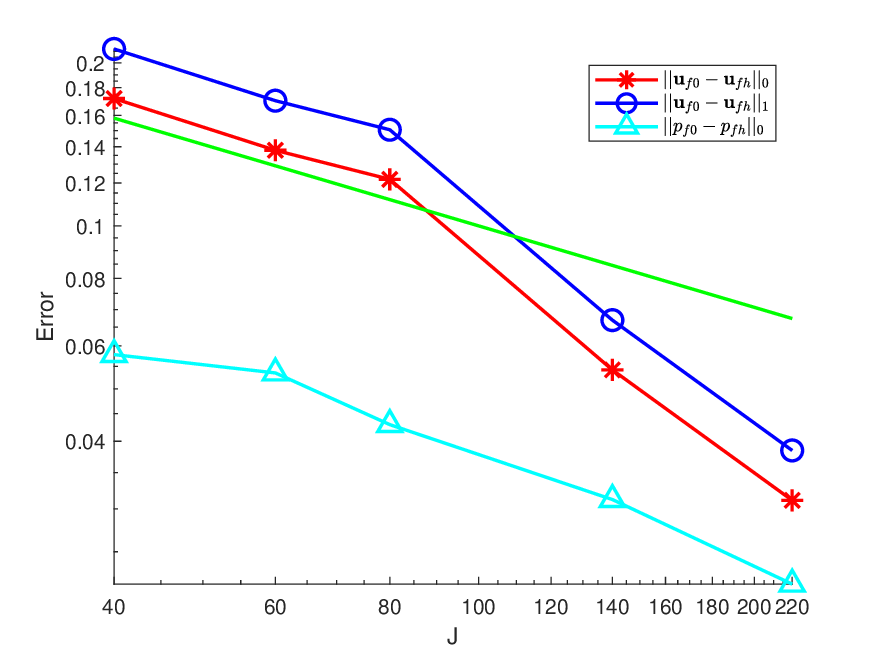} \hspace{5mm}
		\includegraphics[width=0.3\textwidth]{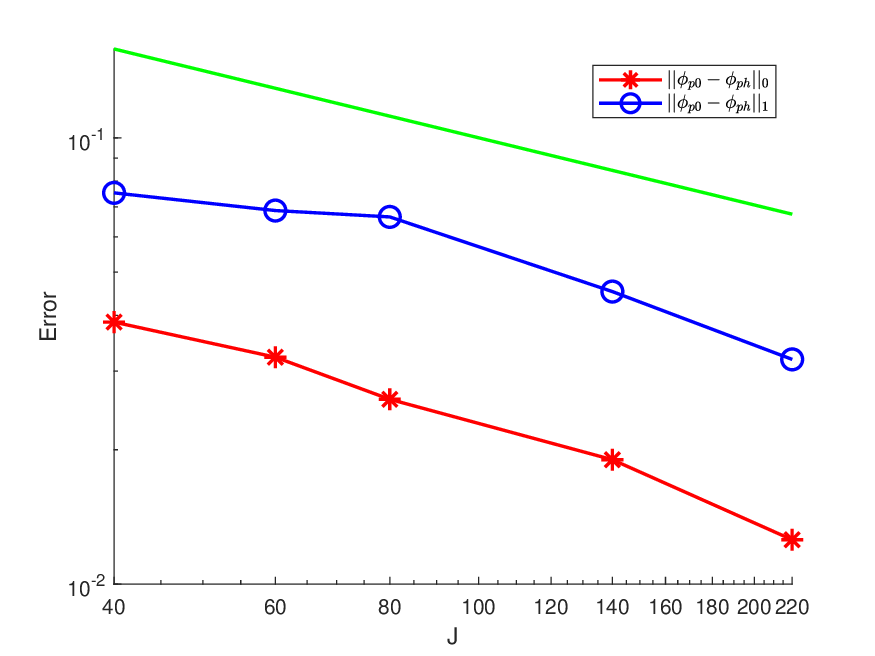}
\caption{The Mento Calo convergence result of the Stokes (left) and Darcy (right) with  $\gamma_{f}^{\ast}$ and $\gamma_{p}^{\ast}$.}
\end{figure} 

\begin{figure}[htbp]
	\centering
		\includegraphics[width=0.27\textwidth]{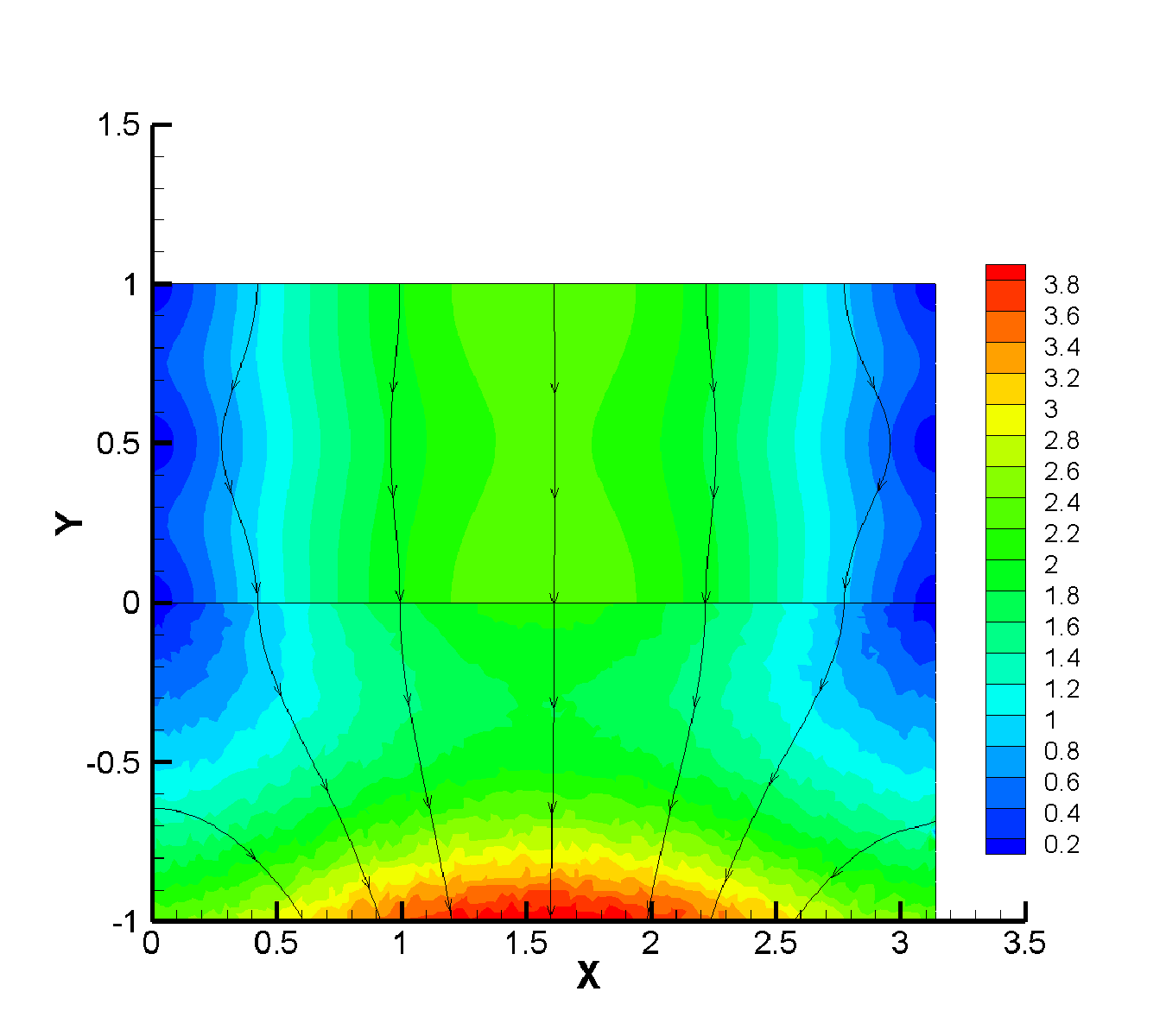}
		\hspace{5mm}
	\includegraphics[width=0.27\textwidth]{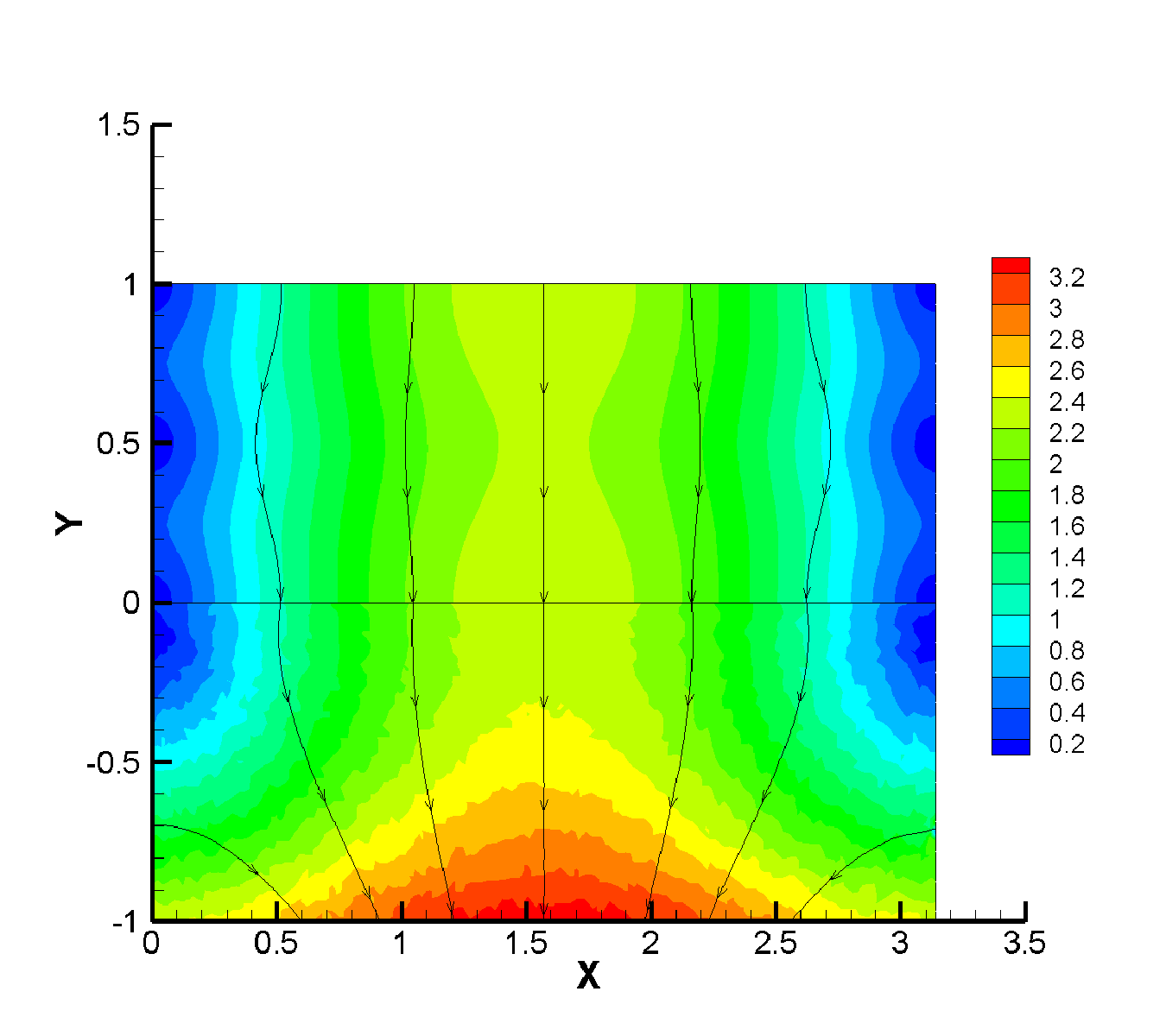}
\hspace{5mm}
\includegraphics[width=0.27\textwidth]{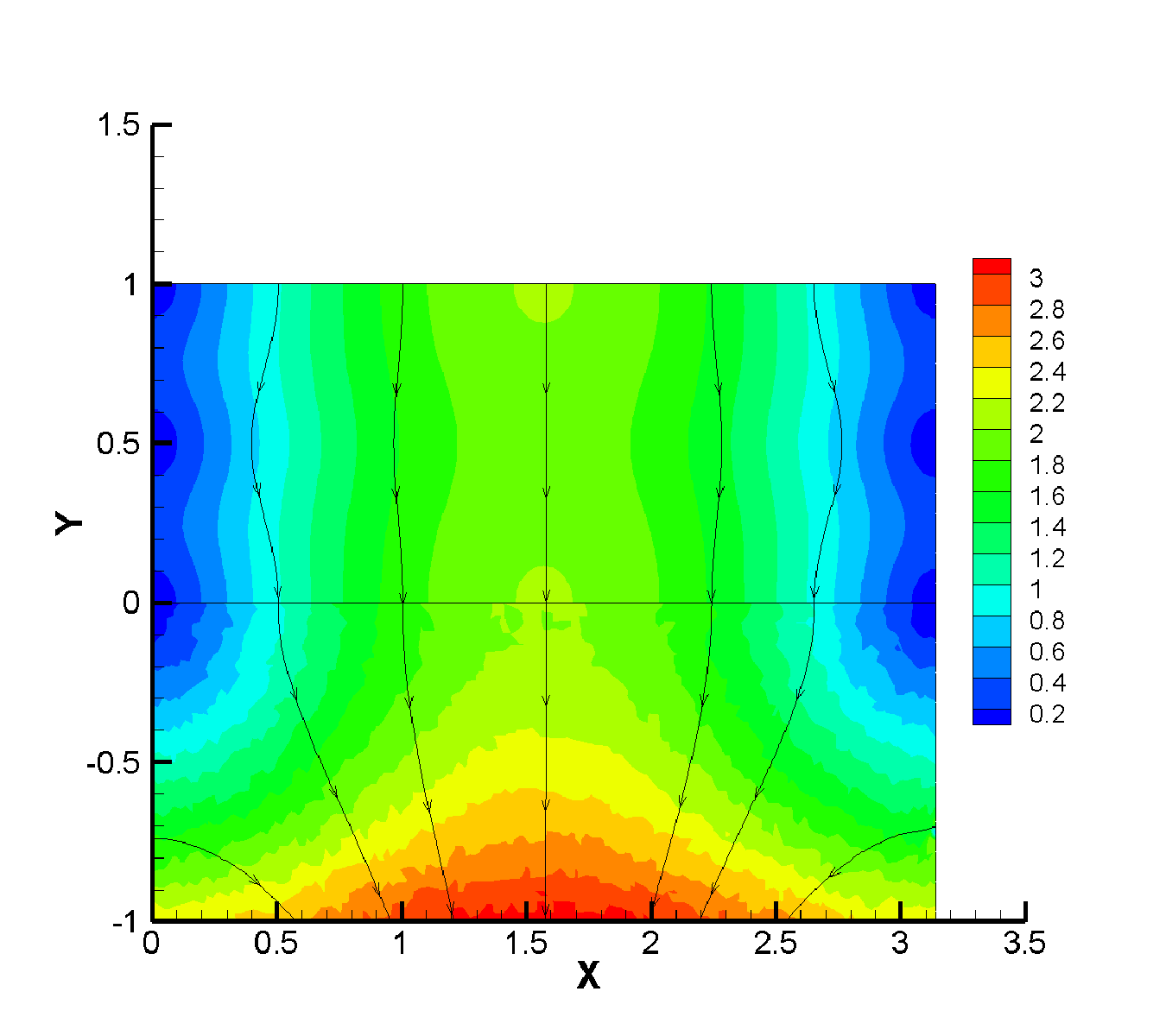}
\caption{The velocity streamlines of velocity expectations of the MC ensemble DDM (left), the traditional DDM (middle) and the MLMC ensemble DDM (right).}\label{SSS}
\end{figure}

We also plot numerical results of the velocity expectations for both domains by the MC ensemble DDM and traditional DDM algorithm. Then, with the choices of $J=80$ realizations and the mesh size $h=\frac{1}{32}$, we present the velocity streamlines in Fig. \ref{SSS} of the random Stokes-Darcy model. Moreover, the third figure of Fig. \ref{SSS} is the velocity result of the MLMC ensemble DDM . We fix $L=2$ and choose the mesh size $h=\frac{1}{2^{3+l}}$ with the corresponding samples number $J=2^{4(L-l)+1}$. In Fig. \ref{SSS}, we can find that the magnitudes of velocity are scattered among the random Stokes-Darcy domain reasonably for both methods. 

	To confirm that MLMC ensemble DDM offers more efficient numerical solutions based on the overall computational effort, we compare its results with those of MC ensemble DDM. We use the ensemble mean of numerical solutions from $J_0=1000$ realizations as our exact solution (expectation). For the MLMC ensemble DDM, we apply the setting $L=2, h=\frac{1}{2^{3+l}}, J=2^{4(L-l)+1}$ as mentioned above. We choose $J=60,80,512$ as the number of MC ensemble DDM samples for a fair comparison, herein, we have calculated the case $J=60, 80$. Table 7.2 displays the errors of both Stokes velocity, pressure and Darcy piezometric head and CPU time. It is clear that the $L_2$ velocity and piezometric head error accuracy of MC ensemble DDM with $512$ samples is similar to that of MLMC ensemble DDM, and MLMC can reduce the CPU time by almost $94.5\%$. The accuracy of the $L_2$ pressure error and the $H^1$ piezometric head error is not very high, but it can still match the error accuracy of MC with 60 samples.
\begin{table}[!h]
	\caption{Comparison of the MC ensemble DDM and MLMC ensemble DDM.}
	\label{MLMCtable}\tabcolsep 0pt \vspace*{-10pt}
	\par
	\begin{center}
		\def\temptablewidth{1.0\textwidth}
		{\rule{\temptablewidth}{1pt}}
		\begin{threeparttable}
			\begin{tabular*}{\temptablewidth}{@{\extracolsep{\fill}}c|c|ccccc|c}
				\hline	
				\hspace{1mm} Method \hspace{1mm} &  $J$ & $||\mathbf u_{f0}-\mathbf u_{fh}||_{0}$ & $||\mathbf u_{f0}-\mathbf u_{fh}||_{1}$ & $||p_{f0}-p_{fh} ||_{0}$   & $||\phi_{p0}-\phi_{ph}||_{0}$ & $||\phi_{p0}-\phi_{ph}||_{1}$ & CPU \\
				\hline
				MLMC  &  512-32-2  & 0.022516     &   0.076512  &        0.055462  &      0.003112  &     0.062549	&	927.71 
				\\
				\hline
				  & 60 &  0.137911   &      0.170296   &     0.053492  &      0.032242  &      0.068710	&		1624.32 \\
				 MC & 80 &  0.121880    &     0.150491  &      0.043936   &     0.025947   &     0.066483	&		2157.66 \\  
				& 512  &    0.021496   &   0.026535  &  0.004455  &    0.002850  &   0.013947  &  16979.3  \\ 
				\hline
			\end{tabular*}%
	\end{threeparttable}
\end{center}
\end{table}

 The convergence effect of the MLMC ensemble DDM method is verified by Table \ref{MLMCtable2} with $L=1,2,3$ corresponding to the mesh size $h=\frac{1}{2^{l+2}}$ and the samples number $J=2^{4(L-l)+1},l=0,\cdots,L$.  We observe that the errors converge to the first order or even higher with respect to $h^L$, which is consistent with Theorem 6.1.
\begin{table}[!h]
	\caption{The convergence effect of MLMC ensemble DDM with different $L$.}
	\label{MLMCtable2}\tabcolsep 0pt \vspace*{-10pt}
	\par
	\begin{center}
		\def\temptablewidth{0.80\textwidth}
		{\rule{\temptablewidth}{0.8pt}}
		\begin{threeparttable}
			\begin{tabular*}{\temptablewidth}{@{\extracolsep{\fill}}c||ccc}
				\hline	
				\hspace{2mm} $L$ \hspace{2mm} &  $||\mathbf u_{f0}-\mathbf u_{fh}||_{1}$  & $||p_{f0}-p_{fh} ||_{0}$   &  $||\phi_{p0}-\phi_{ph}||_{1}$ \hspace{2mm}  \\
				\hline
				\hspace{2mm} 1 \hspace{2mm} & 1.001750    & 0.978971    &  0.392172  \hspace{2mm}
				\\
				\hline
				\hspace{2mm} 2 \hspace{2mm} & 0.519435   &  0.344685    &   0.187737 \hspace{2mm}  \\ 
				\hline
				\hspace{2mm} 3 \hspace{2mm} & 0.073175    & 0.057811  &   0.061237   \hspace{2mm} \\ 
				\hline
			\end{tabular*}%
			\end{threeparttable}
	\end{center}
\end{table}

Finally, 
we compare the computational efficiency of our proposed algorithm and traditional DDM under the selected $J=1,10,20,40,80,160$ realizations and record the respective elapsed CPU times. Let $k_{11}(\vec{x}, \omega), k_{22}(\vec{x}, \omega)$ be the random numbers between 1 and 2. A comparison of the elapsed CPU time of both methods is shown in Table \ref{CPUtime}, which can be clearly that our ensemble DDM algorithm is much faster than the traditional DDM approach.
\begin{table}[!h]
	\caption{The comparison of the elapsed CPU time while the mesh size $h=\frac{1}{64}$.}
	\label{CPUtime}\tabcolsep 0pt \vspace*{-10pt}
	\par
	\begin{center}
		\def\temptablewidth{1.0\textwidth}
		{\rule{\temptablewidth}{1pt}}
		\begin{tabular*}{\temptablewidth}{@{\extracolsep{\fill}}ccccccc}
			\hline
			J&1 & 10 & 20  & 40   & 80 &160\\
			\hline
			traditional DDM    &    0.448   &     4.285   &   8.655  &    18.885  &    36.6  &  69.284
			\\ 
			ensemble DDM   &     0.439  &     2.349   &   4.634 &    10.4  &     20.218&  37.159
			\\
			\hline
		\end{tabular*}%
	\end{center}
\end{table}
%

\subsection{`` 3D Shallow Water'' System with Random Hydraulic Conductivity}
\begin{figure}[htbp]
	\centering
	\subfigure[]{
		\centering\label{3Dshallow}
		\includegraphics[width=45mm,height=55mm]{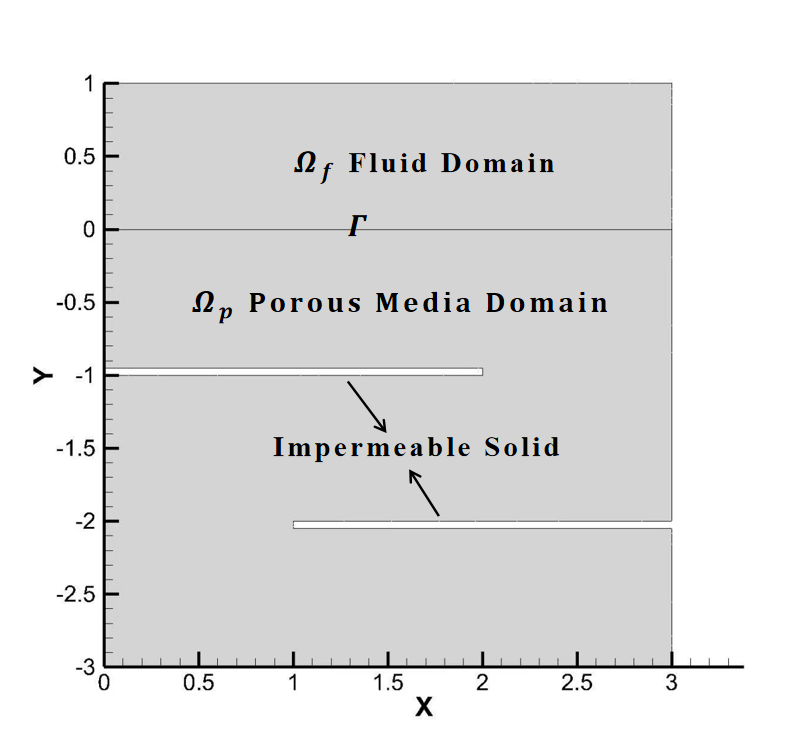}
	}\hspace{5mm}\subfigure[]{
		\centering\label{3Dpp}
		\includegraphics[width=40mm,height=50mm]{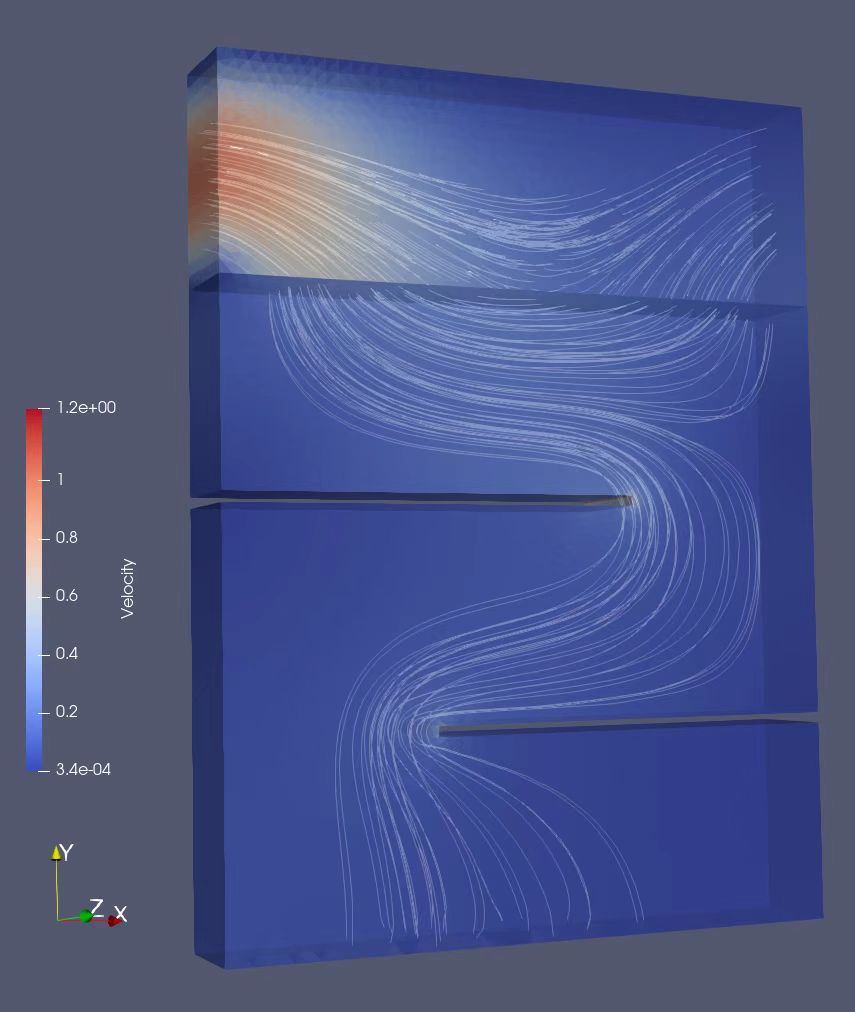}
	}\hspace{5mm}\subfigure[]{
\centering\label{3Dppp}
\includegraphics[width=40mm,height=50mm]{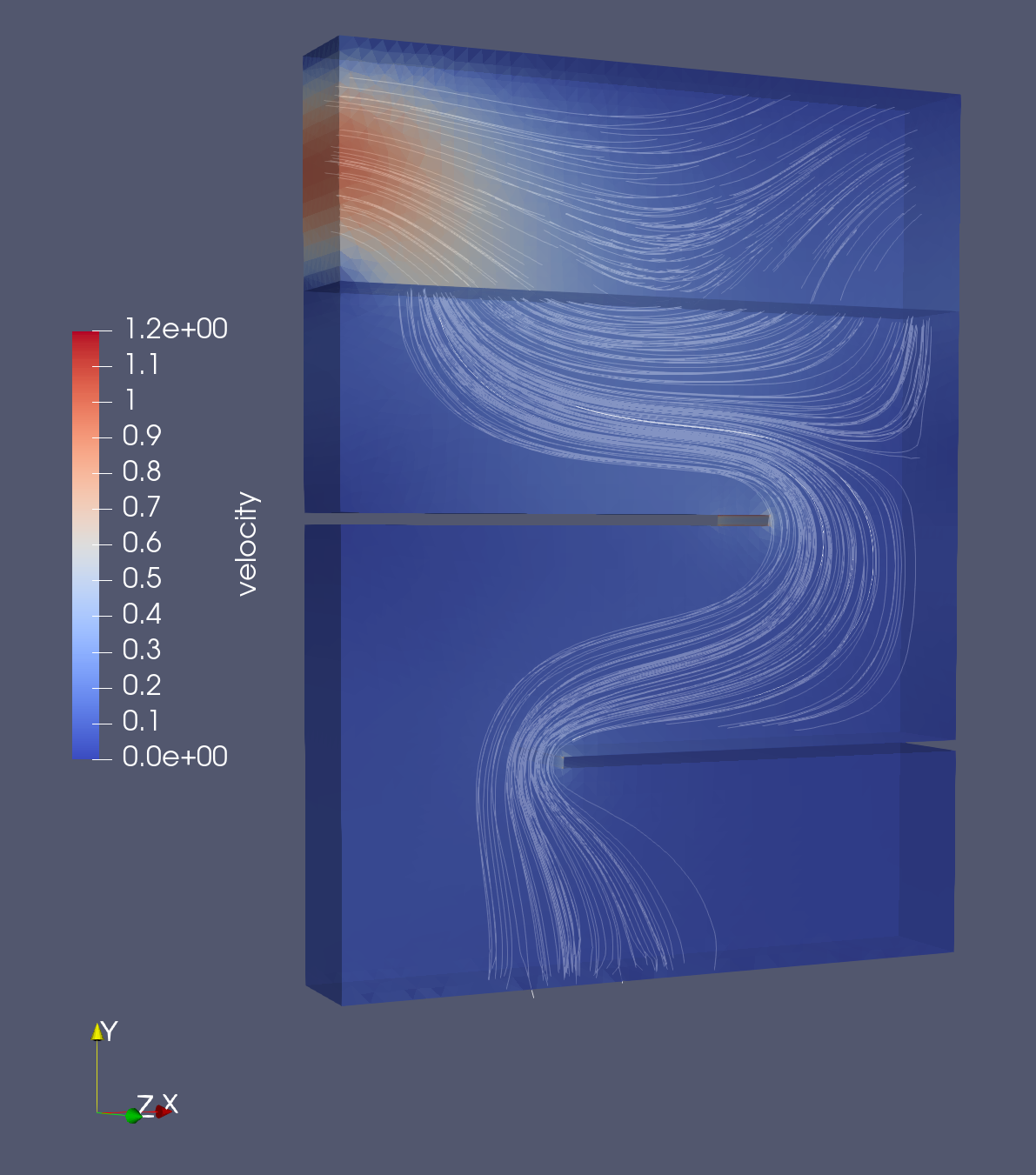}
}
	\caption{(a) The cross-section of the simplified coupling 3D Shallow Water system $\Omega_{f}$ with the porous medium $\Omega_{p}$. (b) The 3D results of the simplified coupling 3D Shallow Water system with the porous medium by MC ensemble DDM. (c) The 3D results of the simplified coupling 3D Shallow Water system with the porous medium by MLMC ensemble DDM (L=2).}
\end{figure}

This experiment is to show the efficiency of the proposed ensemble DDM by simulating a more complicated coupled fluid flow model, known as the "3D Shallow Water" system with the porous medium. We set two impermeable solids for better visualization effects, as illustrated in Fig. \ref{3Dshallow}. Specifically, more details can refer to \cite{Sun05}.
We impose $\mathbf{u}_f = [4y(1-y),0,0]$ as the inflow surface velocity and non-reflective boundary condition on the outlet surface. The rest of the boundaries except the interface of the Stokes subdomain are treated with no-slip boundary conditions. Meanwhile, two impermeable solids of the porous media are constructed with the hydraulic conductivity $k$ of 0.0. Homogeneous Neumann boundary condition $\triangledown\phi_{p}\cdot \mathbf{n}_{p}$ is considered on the vertical ($Y$ direction) face of the porous medium. Homogeneous Dirichlet boundary condition $\phi_{p} = 0$ is applied on the bottom surface of the porous medium. The physical parameter values are taken as $\nu = 1.0, z=0, g=1.0, \mathbf{f}=0$. For the spatial discretization of the free flow region or porous media, we can extrude the 2D(X-Y surface) triangular meshes along the Z direction to form the 3D prism meshes.
We construct the random $\mathbb{K}$ which varies in the vertical direction as follows
\begin{eqnarray*}
	\begin{aligned}
		&\mathbb{K}(\mathbf{x}, \omega)=\left[\begin{array}{ccc}
			k_{11}(\mathbf{x}, \omega) & 0 & 0 \\
			0 & k_{22}(\mathbf{x}, \omega) & 0 \\
			 0 & 0 & k_{33}(\mathbf{x}, \omega)
		\end{array}\right], \quad \text { and } \\
		&k_{11}(\mathbf{x}, \omega)=k_{22}(\mathbf{x}, \omega)=k_{33}(\mathbf{x}, \omega) \\
		&\hspace{13mm}=a_{0}+\sigma \sqrt{\lambda_{0}} Y_{0}(\omega)+\sum_{i=1}^{n_{f}} \sigma \sqrt{\lambda_{i}}\left[Y_{i}(\omega) \cos (i \pi y)+Y_{n_{f}+i}(\omega) \sin (i \pi y)\right],
	\end{aligned}
\end{eqnarray*}
where $\mathbf{x}=(x, y, z)^{T}, \lambda_{0}=\frac{\sqrt{\pi L_{c}}}{2}, \lambda_{i}=\sqrt{\pi} L_{c} e^{-\frac{\left(i \pi L_{c}\right)^{2}}{4}}$ for $i=1, \ldots, n_{f}$ and $Y_{0}, \ldots, Y_{2 n_{f}}$ are uncorrelated random variables with unit variance and zero means. In the following numerical test, we can choose the desired physical correlation length $L_{c}=0.25$ for the random field and $a_{0}=1, \sigma=0.15, n_{f}=3$. The random variables $Y_{0}, \ldots, Y_{2 n_{f}}$ are assumed independent and uniformly distributed in the interval $[-\sqrt{3}, \sqrt{3}] $. With the above assumptions, we note that the corresponding random hydraulic conductivity tensor $\mathbb{K}$ is symmetric positive definite since the random functions $k_{11}(\vec{x}, \omega)$, $k_{22}(\vec{x}, \omega)$ and $k_{33}(\vec{x}, \omega)$ are positive.

We present the velocity vector field and streamline for the MC ensemble DDM with the optimal Robin parameters which are $\gamma_{f}^{\ast}$ and $\gamma_{p}^{\ast}$ respectively, in Fig. \ref{3Dpp}. Meanwhile, the MLMC ensemble DDM with $L=2, h=\frac{1}{2^{2+l}}, J=2^{4(L-l)+1}$ is presented in Fig. \ref{3Dppp}. As shown in Fig. \ref{3Dpp} and Fig. \ref{3Dppp}, the magnitudes of the velocity profile and streamline appear reasonable and present the good performance of the proposed scheme for this 3D model.

\begin{figure}[htbp]
	\centering
	\subfigure[h=1/8]{
		\centering
		\includegraphics[width=40mm,height=30mm]{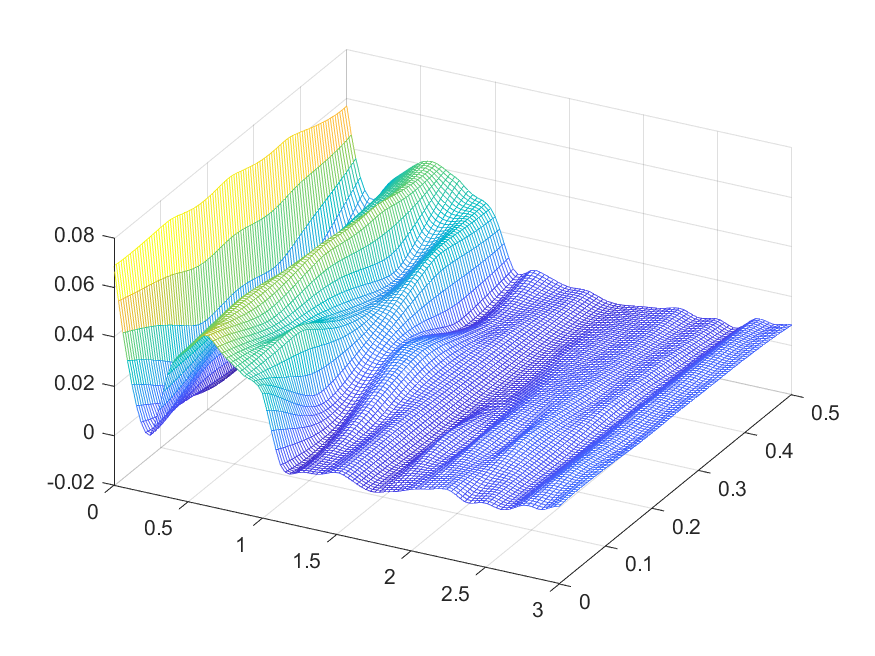}
	}\hspace{3mm}\subfigure[h=1/16]{
		\centering
		\includegraphics[width=40mm,height=30mm]{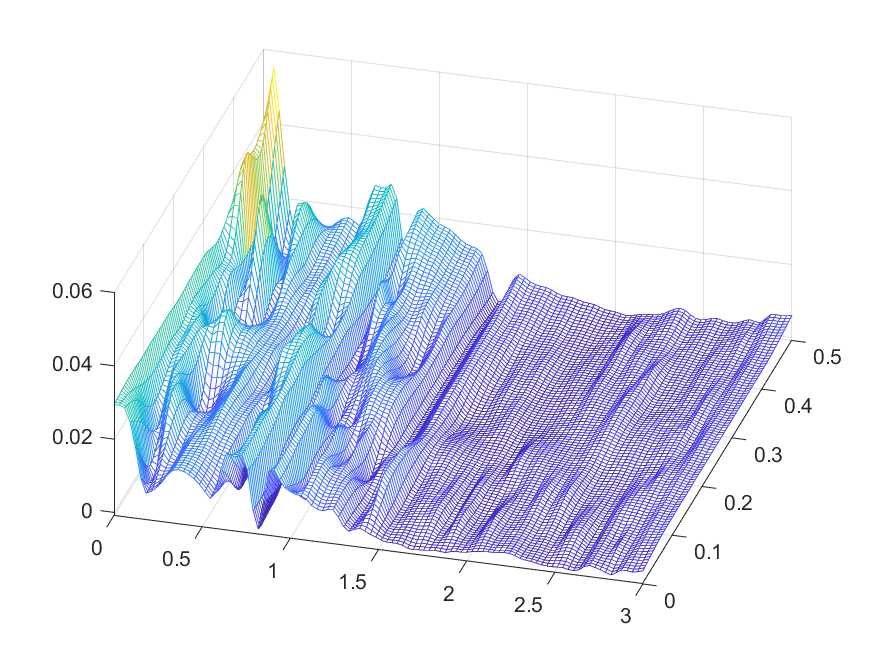}
	}\hspace{3mm}\subfigure[MLMC(L=2)]{
	\centering
	\includegraphics[width=40mm,height=30mm]{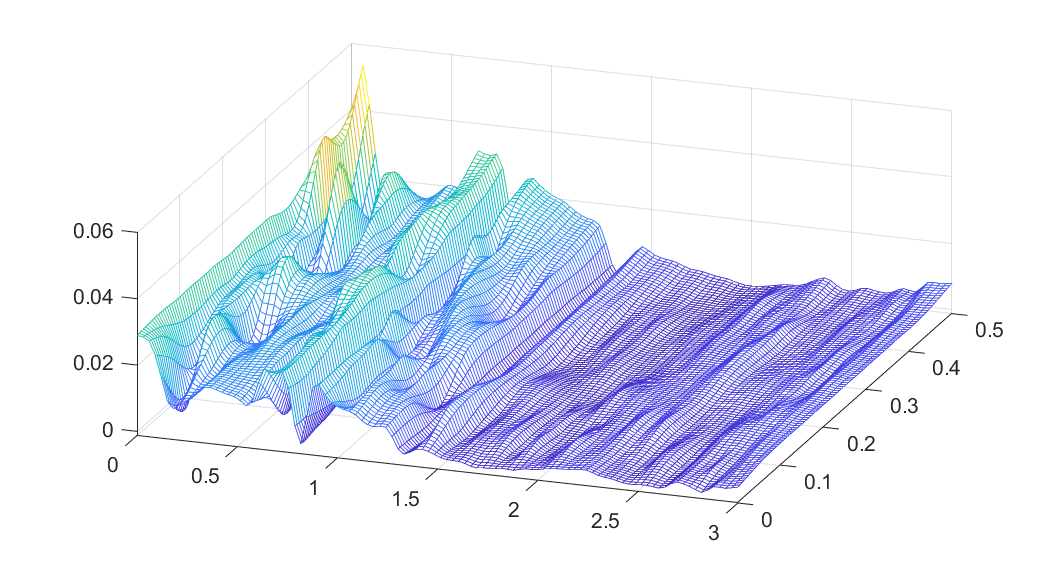}
}
	\caption{(a), (b):The numerical errors on the interface cross-section for different mesh size $h$ by the MC ensemble DDM; (c):The numerical errors on the interface cross-section by the MLMC ensemble DDM  (L=2).}\label{3D3011}
\end{figure}

Furthermore, we check the conservation of mass along the interface on the cross-section of the numerical solution $Y=0.0$ by plotting the errors of $|\mathbf{u}_f\cdot\mathbf{n}_{fh} - \mathbb{K}\triangledown\phi_{ph}\cdot \mathbf{n}_{p}|$ for different mesh scales $h$ in Fig. \ref{3D3011}. By observation, we can find that the continuity of the normal velocity and the approximate error becomes smaller when the mesh turns finer. And the approximate error of the MLMC ensemble DDM has similar accuracy but requires a smaller sample number. 

Through the above experiments, the efficiency of the proposed ensemble DDM algorithm is verified and the property of the optimal Robin parameter is shown. Additionally, the ensemble DDM based on multi-level Monte Carlo exhibits better computational efficiency and accuracy. Moreover, the significant CPU performance, computational efficiency, and complicated flow characteristics are strictly illustrated.

\section{Conclusions}


In this paper, we propose an ensemble Robin-type domain decomposition algorithm for fast-solving random steady-state Stokes-Darcy systems with BJS interface conditions. Our approach involves constructing a single coefficient matrix for all realizations, directly reducing the computational cost by transforming complex and bulky random parameter PDEs into linear systems. Additionally, by employing the Robin-type domain decomposition method, the coupled Stokes-Darcy model is decomposed into Stokes and Darcy subproblems thus the computational consumption can be reduced by realizing parallel computing. We prove the convergence of the algorithm in the continuous case when the Robin parameters satisfy $\gamma_{f}\leq\gamma_{p}$, and the mesh-independent convergence in the case $\gamma_{f}<\gamma_{p}$. Optimized Robin parameters are derived and analyzed to accelerate the convergence of the proposed algorithm. We also discuss finite element approximations for the ensemble domain decomposition algorithm and show both the mesh-dependent and mesh-independent convergence rates when the Robin parameters satisfy $\gamma_{f}\leq\gamma_{p}$. We further propose the MLMC ensemble DDM to significantly lower the computational cost than MC in the probability space, as the number of samples drops quickly when the mesh becomes finer. Finally, we present three numerical experiments to illustrate the efficiency and effectiveness of our ensemble DDM algorithm.

\section{Acknowledgments}
This work was supported by Science and Technology Commission of Shanghai Municipality (Grant Nos.22JC1400900, 21JC1402500, 22DZ2229014).

\end{document}